# Contraction Dynamics in Heterogeneous Spatial Environments

Carlos Barajas, Jean-Jacques Slotine, and Domitilla Del Vecchio *

Department of Mechanical Engineering, Massachusetts Institute of Technology, Cambridge, Massachusetts

**Abstract**

Understanding the asymptotic behavior of reaction-diffusion (RD) systems is crucial for modeling processes ranging from species coexistence in ecology to biochemical interactions within cells. In this work, we analyze RD systems in which diffusion is modeled using the $\theta$-diffusion framework, while the reaction dynamics are spatially varying. We demonstrate that spatial heterogeneity affects the asymptotic behavior of such systems. Using contraction theory, we derive conditions that guarantee the exponential convergence of system trajectories, regardless of initial conditions. These conditions explicitly account for the influence of spatial heterogeneity in both the diffusion and reaction terms. As an application, we study a biochemical system and derive the quasi-steady-state (QSS) approximation, illustrating how spatial heterogeneity modulates the effective binding rates of biomolecular species.

## 1 Introduction

Reaction-diffusion (RD) systems play a pivotal role in modeling biological phenomena, from ecological population dynamics to intracellular biochemical interactions. These systems effectively capture how spatial heterogeneity influences system dynamics [1, 2]. The asymptotic behavior of RD systems is fundamental in understanding species coexistence, disease persistence, environmental sustainability, and in designing genetic circuits for controlling biological processes [3–7]. When diffusion occurs on a much faster timescale than reaction dynamics, model reduction techniques indicate that solutions to the partial differential equation (PDE) system can be approximated by ordinary differential equations (ODEs) after a rapid transient phase [6, 8, 9]. Such approximations simplify analysis and facilitate design.

However, many existing studies focus on specific conditions, typically assuming spatially homogeneous reaction dynamics and Fickian diffusion, where species migrate opposite to concentration gradients. These studies aim to derive conditions ensuring spatial homogenization of solutions [10–13], providing rigorous support for the well-mixed assumption commonly used in modeling biomolecular processes [14]. Yet, accurately modeling diffusion in spatially heterogeneous environments requires more general diffusion models, such as dispersal models [1, 15]. For instance, in bacterial cells, diffusion involves flux terms comprising components that oppose concentration gradients and components that drive species outward from densely packed chromosomal regions [2]. Under rapid diffusion, it has been shown that asymptotic solutions in such systems are not homogeneously distributed and can modulate binding affinities between reactants [6, 16]. Moreover, reaction dynamics can explicitly depend on spatial variables, for instance, the parameters that govern disease spread in epidemiological models [17–19]. Because these complexities, the asymptotic behavior of systems incorporating dispersal diffusion and spatially varying reaction dynamics has not been rigorously analyzed.

In this study, we employ nonlinear contraction theory [20, 21] to investigate the asymptotic behavior of reaction-diffusion (RD) systems with $\theta$-diffusion— an extensible, one-parameter family of dispersal models [15] that generalizes Fickian diffusion— and spatially varying reaction dynamics. We derive conditions based on the weighted spatial average dynamics, the second eigenvalue of the diffusion operator, and the Jacobian of the reaction dynamics to ensure that all solutions of the system converge exponentially fast to each other, independent of initial conditions. These conditions explicitly reveal how spatial heterogeneity in both diffusion and reaction terms influences the system's asymptotic behavior. We apply our theoretical results to a biochemical reaction system, demonstrating the quasi-steady-state (QSS) approximation and highlighting how spatial heterogeneity modulates effective binding rates of biomolecular species.

The paper is organized as follows. In Section 3, we present motivating examples that illustrate how spatial heterogeneity can destabilize RD systems and affect bimolecular binding rates. In Section 4, we introduce our main theoretical results, providing conditions that guarantee exponential convergence of solutions to the RD system. In Section 5, we revisit the motivating examples and present additional cases to demonstrate the applicability of our results. Finally, we conclude with remarks summarizing our findings and discussing their implications for modeling and analysis of RD systems.

---

*Correspondence: ddv@mit.edu

## 2 Notation

In this work, $t \in [0, \infty)$ and $x \in \Omega$ denote time and space, respectively. Let $\Omega \subset \mathbb{R}^k$ be a bounded, connected domain with smooth boundary $\partial \Omega$ and closure $\bar{\Omega}$. For simplicity, we assume that the volume of $\Omega$ is normalized to 1, i.e., $|\Omega| := \int_\Omega dx = 1$. We denote the $(j,k)$-th element of a matrix $A$ by $A^{j,k}$.

We introduce the $L^2$-norm on the space $L^2(\Omega, \mathbb{R}^n)$, which applies to functions mapping $\Omega$ to $\mathbb{R}^n$, as

$$\|u(t,\cdot)\| := \left( \int_\Omega u^T(t,x) u(t,x) dx \right)^{1/2},$$

where $u^T(t,x)$ denotes the transpose of the vector $u(t,x) \in \mathbb{R}^n$. By the fact that $|\Omega| = 1$, when $u(t)$ is independent of $x$, this norm reduces exactly to the standard Euclidean norm in $\mathbb{R}^n$. In this work, if a vector quantity is set to zero, it implies that the condition holds element-wise for all components of the vector. For a vector $v \in \mathbb{R}^n$, we denote by $\text{diag}(v)$ the $n \times n$ diagonal matrix whose diagonal entries are $v_1, v_2, \ldots, v_n$ and whose off-diagonal entries are all zero. Finally, we let $I_{n \times n}$ denote the identity matrix of dimension $n \times n$. We write $u(t,x) \to u^*(t,x)$ *exponentially* if there exists a finite constants $C > 0$ and $\lambda > 0$ such that

$$\|u(t,x) - u^*(t,x)\| \leq C e^{-\lambda t}, \quad \forall t \geq 0.$$

For this work, a function is considered smooth if every required derivative (or partial derivative) exists and is continuous.

## 3 Motivation

This work is motivated by the observation that spatial heterogeneity can significantly impact the asymptotic behavior of RD systems. In particular, we illustrate how spatial variations in both the diffusion environment and reaction terms can destabilize these systems. Additionally, we demonstrate how spatial heterogeneity modulates the binding rates of biomolecular species within individual cells.

### 3.1 Spatially varying reaction dynamics

Consider a scalar reaction-diffusion system with Fickian diffusion, where the linear reaction dynamics exhibit spatial variation. The system evolves under Neumann boundary conditions and is described by:

$$\frac{\partial z(t,x)}{\partial t} = \frac{\epsilon}{\pi^2} \frac{\partial^2 z(t,x)}{\partial x^2} + a(x) z(t,x), \quad x \in (0,1), \qquad (1)$$

$$\frac{\partial z(t,x)}{\partial x} = 0, \quad x \in \{0,1\},$$

where $\epsilon > 0$. Here, $a(x) = -\epsilon + \sin(\omega x) - \int_0^1 \sin(\omega x') \, dx'$ represents a spatially varying reaction term, and $\bar{a} := \int_0^1 a(x) \, dx = -\epsilon$. As $\omega$ increases, the heterogeneity of $a(x)$ across the domain $x \in (0,1)$ becomes more pronounced while it's space average, $\bar{a}$, remains constant.

The time evolution of $\|z\|$ for various values of $\omega$ is depicted in Fig. 1(a). For small $\omega$, $\|z\|$ decays exponentially, indicating stability, whereas for larger $\omega$, $\|z\|$ grows exponentially, signifying instability. This demonstrates how increasing spatial heterogeneity in the reaction term can induce instability in the system dynamics. The relationship between $\omega$ and the slope of the logarithm of $\|z\|$ after an initial transient ($t > 80$) is shown in Fig. 1(b). We observe stability for $\omega \lesssim 0.3$. In this work, we analyze the role of spatial heterogeneity in the reaction dynamics in the asymptotic behavior of RD systems.

## 3.2 Spatially heterogeneous diffusion dynamics

Consider a two-state reaction-diffusion system where $z(t,x) = [z_1(t,x),\ z_2(t,x)]^\top$ represents the state vector, with linear, spatially independent reaction dynamics in a one-dimensional spatial domain $\Omega = (0,1)$. The system is governed by:

$$\frac{\partial z(t,x)}{\partial t} = -\frac{d}{dx}\left[J(x,z)\right] + Az(t,x), \tag{2}$$
$$J(x,z) = 0, \quad x \in \{0,1\},$$

where $J(x,z)$ denotes the flux vector, and $A$ is the reaction matrix. The reaction matrix $A$ is given by:

$$A = \begin{bmatrix} -1 & 1 \\ -1 & \frac{1}{2} \end{bmatrix}.$$

Let $D = \frac{\zeta}{\pi^2}\left[10,\ \frac{1}{4}\right]^T$, where $\zeta > 0$ controls the overall speed of diffusion. We define flux vector $J(x,z) = [J^{1,1}(x,z),\ J^{2,1}(x,z)]^\top$ as follows:

1. The first component, $J^{1,1}(x,z)$, obeys Fick's law of diffusion:

$$J^{1,1}(x,z) = -D^{1,1}\frac{\partial z_1(t,x)}{\partial x}. \tag{3}$$

2. The second component, $J^{2,1}(x,z)$, captures diffusion in a heterogeneous environment and is given by:

$$J^{2,1}(x,z) = -D^{2,1}\left[v_r(x)\frac{\partial z_2(t,x)}{\partial x} - z_2(t,x)\frac{\partial v_r(x)}{\partial x}\right], \tag{4}$$

where $v_r(x) \in (0,1]$ for all $x \in [0,1]$ represents spatial heterogeneity in the environment [2, 15]. The function $v_r(x)$ models the available volume fraction and depends on a parameter $r$ that characterizes the degree of heterogeneity. It is plotted for several values of $r$ in Fig. 2(a). Details on the physical interpretation of $v_r(x)$ in the context of intracellular biochemical reactions is discussed in SI section 7.1.

The flux term (4) includes a component in the same direction as Fickian diffusion (opposite to $\frac{\partial z_2(t,x)}{\partial x}$) and an additional component in the direction of $\frac{\partial v_r(x)}{\partial x}$, accounting for spatial variation in the environment. Note that when $r = 0$, we have $v_r(x) = 1$ (Fig. 2(a)), and the flux (4) reduces to Fick's diffusion (3).

Let $\bar{z}(t) := \int_0^1 z(t,x)\,dx$ denote the spatial average of $z(t,x)$, and define:

$$z^*(t,x) := \mathrm{diag}\left(\begin{bmatrix} 1, & \hat{v}_r(x) \end{bmatrix}^T\right)\bar{z}(t) \tag{5}$$

where $\hat{v}_r(x) := \dfrac{v_r(x)}{\int_0^1 v_r(x)\,dx}$. The function $z^*(t,x)$ represents the no-flux solution satisfying $J(x, z^*(t,x)) = 0$.

Define the fluctuations from the no-flux solution as $z^\perp(t,x) := z(t,x) - z^*(t,x)$. By integrating both sides of (2) over the spatial variable and applying the boundary conditions, one can show that:

$$\frac{d\bar{z}(t)}{dt} = A\bar{z}(t), \tag{6}$$
$$\frac{\partial z^\perp(t,x)}{\partial t} = -\frac{d}{dx}\left[J(x,z^\perp)\right] + Az^\perp(t,x) + A\left(\mathrm{diag}\left(\begin{bmatrix} 1, & \hat{v}_r(x) \end{bmatrix}^T\right) - I_{2,2}\right)\bar{z}(t), \quad x \in (0,1),$$
$$J(x,z^\perp) = 0, \quad x \in \{0,1\}.$$

Notice that the dynamics of $\bar{z}(t)$ are independent of $z^\perp(t,x)$, and since $A$ is Hurwitz (its eigenvalues have negative real parts), it implies that $\bar{z}(t) \to 0$. Conversely, the dynamics of $z^\perp(t,x)$ generally depend on $\bar{z}(t)$, except when $v_r(x) = 1$ for all $x \in (0,1)$, i.e., when there is no spatial heterogeneity in the diffusion dynamics. As we will demonstrate in this work, for general nonlinear systems with spatial heterogeneity in the diffusion and reaction dynamics, the space-averaged state and the fluctuations from the no-flux solution are coupled. This coupling necessitates more stringent conditions than those previously established for homogeneous systems [10].

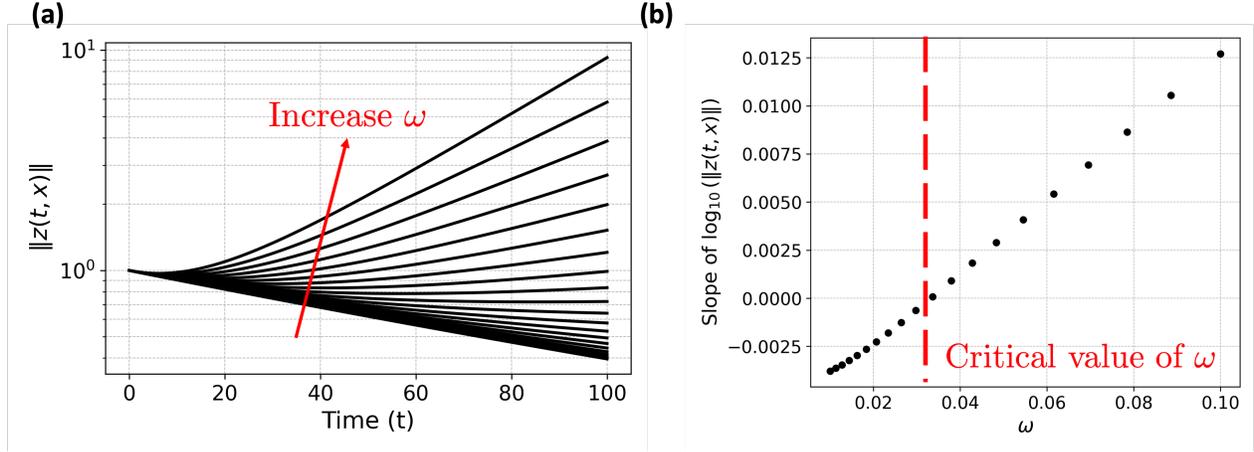

Figure 1: (a) Time evolution of $\|z(t,x)\|$ for various values of $\omega$ in the system described by (1) when $\epsilon = 1 \times 10^{-2}$. (b) The relationship between $\omega$ and the slope of $\log \|z(t,x)\|$ after an initial transient ($t > 80$), showing stability for $\omega \lesssim 0.3$. For full simulation detail see SI Section 7.2.

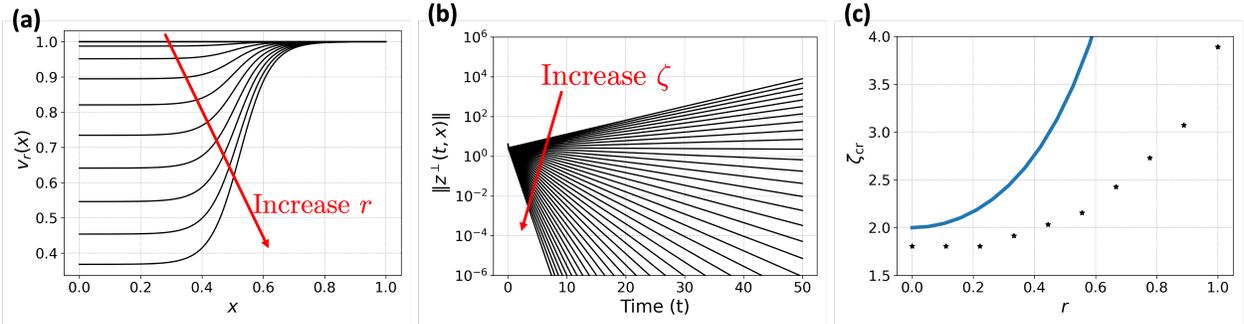

Figure 2: (a) The available volume profile $v_r(x)$ appearing in the flux term (4), as a function of $x$ for several values of $r$. (b) Time evolution of $\|z^\perp(t,x)\|$ for various values of $\zeta$ when $r = 0$ in (4) in the system (2). (c) The critical value of $\zeta$, $\zeta^*(r)$, determined through simulations (black line) and the theoretical upper-bound approximation (blue line) given by (43). For full simulation detail see SI Section 7.3.

The temporal evolution of $\|z^\perp(t,x)\|$ is shown in Fig. 2(b) for several values of $\zeta$ when $r = 0$ in (4), starting from the same initial conditions. It is observed that for small values of $\zeta$ (slow diffusion), $\|z^\perp(t,x)\| \to \infty$, indicating instability. For large values of $\zeta$ (fast diffusion), we observe that $\|z^\perp(t,x)\| \to 0$, indicating stability. This implies that slow diffusion can destabilize $z(t,x)$, while fast diffusion stabilizes it—a phenomenon well-known in the study of Turing patterns [22, 23].

The critical value of $\zeta$, denoted as $\zeta_{\text{cr}}(r)$, at which the dynamics transition from unstable to stable, increases with $r$ (Fig. 2(c)). This implies that greater spatial heterogeneity (larger $r$) requires faster diffusion to stabilize the $z^\perp$ dynamics.

Furthermore, when $z^\perp(t,x) \to 0$, it implies convergence to the no-flux solution $z^*(t,x)$ given by (5) exponentially fast, highlighting a more general phenomenon than spatial homogenization. In this work, we will rigorously establish conditions under which such convergence occurs and analyze its implications for the overall system dynamics.

### 3.3 Spatial heterogeneity modulates binding rates in biomolecular reactions

Consider the intracellular process of mRNA translation, where mRNA (m) binds to a ribosome (R) to form an mRNA-ribosome complex ($c_r$), which then initiates protein (P) synthesis. This process is modeled by the following biochemical reaction scheme:

$$\text{m} + \text{R} \underset{d}{\overset{a}{\rightleftharpoons}} \text{c}_r \overset{\kappa}{\longrightarrow} \text{P} + \text{m} + \text{R} \tag{7}$$

where $a$ and $d$ are the association and dissociation rate constants, respectively, for ribosome-mRNA binding, and $\kappa$ is the catalytic rate constant for protein synthesis. Let $m(t,x)$, $R(t,x)$, and $c_r(t,x)$ denote the concentrations of mRNA, ribosomes, and mRNA-ribosome complexes, respectively, at time $t$ and spatial position $x$ inside the cell. The spatial variable $x$ takes values in the general spatial domain $\Omega \subset \mathbb{R}^3$, normalized such that $|\Omega| := \int_\Omega dx = 1$, representing the volume of the cell.

The biomolecular species described in equation (7) are assumed to diffuse freely within the cell. Their flux dynamics are modeled by spatially heterogeneous equations, which are high-dimensional analogs of equation (4). Specifically, the fluxes are given by:

$$\begin{aligned}
J_m &= D_m\big[-v_{r_m}(x)\nabla m(t,x) + m(t,x)\nabla v_{r_m}(x)\big], \\
J_r &= d\big[-v_{r_r}(x)\nabla R(t,x) + R(t,x)\nabla v_{r_r}(x)\big], \\
J_c &= D_c\big[-v_{r_c}(x)\nabla c_r(t,x) + c_r(t,x)\nabla v_{r_c}(x)\big],
\end{aligned} \tag{8}$$

where $D_i$ represents the diffusion coefficient, and $v_i(x)$ is the spatially varying available volume profile across $\Omega$. Here, the subscript $i$ corresponds to $r_m$, $r_r$, or $r_c$, which denote the radii of gyration of mRNA, ribosomes, and the mRNA–ribosome complex, respectively.

These equations describe the fluxes of biomolecular species within the three-dimensional domain $\Omega$, accounting for spatial heterogeneity in the available volume profiles. Each flux includes a component driven by concentration gradients and another that moves toward regions of higher available volume. The available volume profiles, $v_i(x)$, represent the fraction of the cell accessible for diffusion, determined by the spatially heterogeneous structure of the chromosomal DNA mesh, which restricts diffusion [2]. As shown qualitatively in Fig. 2(a), $v_i(x)$ depends on the radius of gyration $r$ of the diffusing biomolecule, capturing its interaction with the DNA mesh (see SI Section 7.1 for further details).

Let overbars denote the spatially averaged concentrations (e.g., $\bar{m}(t) := \int_\Omega m(t,x)\,dx$). In [6], it was shown that in the limit where diffusion is arbitrarily faster than the reaction dynamics, the concentrations tend to:

$$\begin{aligned}
m(t,x) &\to \bar{m}(t)\hat{v}_{r_m}(x), \\
R(t,x) &\to \bar{R}(t)\hat{v}_{r_r}(x), \\
c_r(t,x) &\to \bar{c}_r(t)\hat{v}_{r_c}(x),
\end{aligned} \tag{9}$$

where this convergence is exponentially fast. Here, $\hat{v}_i(x) := \dfrac{v_i(x)}{\int_\Omega v_i(x)\,dx}$, with the subscript $i$ corresponds to $r_m$, $r_r$, or $r_c$. In this work, we will show that this approximation holds even when diffusion and reactions occur on the same timescale. Furthermore, we will demonstrate that the ribosome-mRNA complex dynamics approach:

$$\bar{c}_r(t) \to \theta \frac{\bar{m}(t)\,\bar{R}(t)}{K}, \tag{10}$$

where

$$\theta = \frac{\int_\Omega v_{r_m}(x) v_{r_r}(x)\,dx}{\left(\int_\Omega v_{r_m}(x)\,dx\right)\left(\int_\Omega v_{r_r}(x)\,dx\right)}, \quad K = \frac{d+\kappa}{a}. \tag{11}$$

This approximation is known as the quasi-steady-state (QSS) approximation [14], and the quantity $\theta$, known as the Binding Correction Factor (BCF), quantifies how spatial heterogeneity in diffusion affects the effective binding rate between mRNA and ribosomes [6]. When $\theta = 1$, equation (10) coincides with a well-mixed model in a homogeneous environment. The BCF increases when the freely diffusing species are concentrated into specific regions of the cell due to a highly dense chromosomal mesh concentrated near mid-cell that forces them to interact within a reduced volume.

## 4 Results

This section introduces a one-parameter family of diffusive flux models that generalizes Equation (4) to capture diffusion in heterogeneous environments. We define the linear differential operator derived from this diffusive flux and demonstrate some of its spectral properties. We then present the main results of the paper, which establish the contraction of RD systems with diffusion modeled by this linear operator and spatially varying reaction dynamics. We illustrate how classic tools from contraction theory, such as virtual dynamics, hierarchal/feedback combinations, and partial contraction, are employed.

## 4.1 Diffusion Operator to Model Heterogeneous Environments

For $\theta \in [0,1]$, let $d : \Omega \to \mathbb{R}_+$ be a positive, smooth function representing a diffusivity coefficient, and let $y : \Omega \to \mathbb{R}$ represent a concentration or state variable. Consider the flux-diffusion model defined by

$$J(\theta, d(x), y(x)) := -d(x)\nabla y(x) + (2\theta - 1)y(x)\nabla d(x) = -d^{2\theta}(x)\nabla\left[d^{1-2\theta}(x)y(x)\right]. \tag{12}$$

We define the associated differential operator

$$L(\theta, d, y) := -\operatorname{div}\left(J(\theta, d(x), y(x))\right), \quad \forall x \in \Omega, \tag{13}$$

subject to the no-flux boundary condition

$$J(\theta, d(x), y(x)) \cdot \nu = 0, \quad \forall x \in \partial\Omega, \tag{14}$$

where $\nu$ denotes the outward unit normal vector on $\partial\Omega$.

The null space of the operator $L(\theta, d, y)$ consists of functions of the form $y(x) = c\psi_{\theta,d}(x)$ for an arbitrary constant $c$, where

$$\psi_{\theta,d}(x) = \frac{d^{2\theta-1}(x)}{\int_\Omega d^{2\theta-1}(x)\,dx}. \tag{15}$$

**Lemma 1.** *Suppose $y(x)$ satisfies the no-flux boundary condition (14), and let $\bar{y} := \int_\Omega y(x)\,dx$ represent the average of $y$ over $\Omega$. Let $\psi_{\theta,d}(x)$ be defined by (15), and define the orthogonal component $y^\perp(x) := y(x) - \psi_{\theta,d}(x)\bar{y}$, which represents the deviation of $y$ from its projection onto the null space. Let $L(\theta, d, y)$ be the differential operator defined in (13). Then the following inequality holds:*

$$\frac{\int_\Omega y^\perp(x)\psi_{\theta,d}^{-1}(x)L(\theta, d, y^\perp)\,dx}{\int_\Omega y^\perp(x)\psi_{\theta,d}^{-1}(x)y^\perp(x)\,dx} \leq -\lambda_{\theta,d}, \tag{16}$$

*for some constant $\lambda_{\theta,d} > 0$ that depends on $\Omega$ and properties of $d$. Furthermore, let $\lambda^*$ be the first non-zero Neumann eigenvalue for the Laplacian on the domain $\Omega$. Then,*

$$\lambda_{\theta,d} \geq \lambda^* \frac{\min_{x\in\bar{\Omega}} d^{2\theta}(x)}{\max_{x\in\bar{\Omega}} d^{2\theta-1}(x)}. \tag{17}$$

*Proof.* The proof relies on the self-adjointness of $L(\theta, d, y)$ with respect to the inner product weighted by $\psi_{\theta,d}^{-1}(x)$, that is,

$$\langle u, v \rangle_\psi := \int_\Omega u(x)\psi_{\theta,d}^{-1}(x)v(x)\,dx,$$

and the min-max principle of elliptic operators ([24] equation 1.37). See SI Section 7.4 for full details. □

The diffusion model (12) is known as the $\theta$-diffusion model [15]. The model (12) coincides with Fickian diffusion when $\theta = \frac{1}{2}$ or $d(x)$ is spatially homogeneous, where the flux is in the opposite direction of $\nabla y$. When $\theta > \frac{1}{2}$ (respectively, $\theta < \frac{1}{2}$) and $d(x)$ is spatially heterogenous, there is an additional flux component in the same (respectively, opposite) direction as $\nabla d(x)$. In mathematical ecology, these are known as biased dispersal models, where the fitness of individuals depends heterogeneously on the environment [1].

## 4.2 Certifying Contraction of Spatially Heterogeneous RD Systems

In this section, we introduce the notion of contraction and provide sufficient conditions that guarantee the contraction of dynamical with spatially heterogeneous reaction dynamics and diffusion described by (13) . We then present two corollaries that extend these results to a broader class of systems, relax the conditions, and demonstrate cases where contraction ensures that the deviations from the no-flux solution ($y(x) \propto d^{2\theta-1}(x)$ in (12)) decay exponentially.

**Definition** 1. *A dynamical system is considered* contracting *if any two trajectories, originating from distinct initial conditions, exhibit exponential convergence toward each other. The* contraction rate *denotes a uniform lower bound on the rate at which these trajectories converge.*

Let $\bar{w}(t) \in \mathbb{R}^n$ and $z^\perp(t,x) \in \mathbb{R}^m$, such that $\int_\Omega (z^\perp(t,x) dx = 0$ for all $t \geq 0$, and suppose that the functions $f_1(t,\bar{w})$, $g_1(t,x,z^\perp)$, $f_2(t,x,z^\perp)$, and $g_2(t,x,\bar{w})$ are smooth. For $\Theta \in \mathbb{R}^m$ and $D(x) \in \mathbb{R}^m$ such that $\Theta^i \in [0,1]$, $D^i(x) > 0$, and $D^i(x)$ are smooth for all $i = 1,\ldots, m$, let

$$\Psi_{\Theta,D}(x) = \text{diag}\left(\psi_{\Theta^1, D^1}(x), \ldots, \psi_{\Theta^m, D^m}(x)\right),$$

where each $\psi_{\Theta^i, D^i}(x)$ is given by equation (15).

Consider the system:

$$\begin{aligned}
\frac{d\bar{w}(t)}{dt} &= f_1(t,\bar{w}) + \int_\Omega g_1(t,x',z^\perp) \, d', \quad \forall t > 0, \\
\frac{\partial z^\perp(t,x)}{\partial t} &= \mathcal{L}_\Theta(D, z^\perp) + f_2^\perp(t,x,z^\perp) + g_2^\perp(t,x,\bar{w}), \quad \forall x \in \Omega, \quad \forall t > 0, \\
J\left(\Theta^i, D^i, z^{\perp,i}\right) \cdot \nu &= 0, \quad \forall x \in \partial\Omega, \quad \forall t > 0, \quad \forall i = 1,\ldots, m,
\end{aligned} \quad (18)$$

where

$$\mathcal{L}_\Theta^i(D, z^\perp) := L\left(\Theta^i, D^i(x), z^{\perp,i}\right),$$

as given by equation (13),

$$f_2^\perp(t,x,z^\perp) := f_2(t,x,z^\perp) - \Psi_{\Theta,D}(x) \, \bar{f}_2(t,z^\perp),$$

$$g_2^\perp(t,x,\bar{w}) := g_2(t,x,\bar{w}) - \Psi_{\Theta,D}(x) \, \bar{g}_2(t,\bar{w}),$$

and

$$\bar{f}_2(t,z^\perp) := \int_\Omega f_2(t,x,z^\perp) \, dx, \quad \bar{g}_2(t,\bar{w}) := \int_\Omega g_2(t,x,\bar{w}) \, dx.$$

**Theorem** 1. *Consider the system given by equation* (18), *and let* $\Lambda_{\Theta,D} = \text{diag}(\lambda_{\Theta^1, D^1}, \ldots, \lambda_{\Theta^m, D^m})$, *where each* $\lambda_{\Theta^i, D^i}$ *is defined by equation* (16). *Suppose there exist convex regions* $\bar{\chi}_w \subseteq \mathbb{R}^n$ *and* $\chi^\perp \subseteq \mathbb{R}^m$, *constants* $\lambda_1, \lambda_2, \beta > 0$, *and uniformly positive definite matrices* $M_1(t) \in \mathbb{R}^{n \times n}$ *and* $M_2(t,x) \in \mathbb{R}^{m \times m}$, *where* $M_2(t,x) := \Gamma(t)\Psi_{\Theta,D}^{-1}(x)$ *and* $\Gamma(t)$ *is a diagonal matrix, such that:*

1.
$$\frac{1}{2}\left[\dot{M}_1 + M_1 \frac{\partial f_1}{\partial \bar{w}} + \left(\frac{\partial f_1}{\partial \bar{w}}\right)^\top M_1\right] \leq -\lambda_1 M_1, \quad \forall t \geq 0, \forall \bar{w} \in \bar{\chi}_w. \quad (19)$$

2.
$$\frac{1}{2}\left[\dot{M}_2 + M_2 \left(\frac{\partial f_2}{\partial z^\perp} - \Lambda_{\Theta,D}\right) + \left(\frac{\partial f_2}{\partial z^\perp} - \Lambda_{\Theta,D}\right)^\top M_2\right] \leq -\lambda_2 M_2, \quad \forall t \geq 0, \forall x \in \Omega, \forall z^\perp \in \chi^\perp. \quad (20)$$

3.
$$\int_\Omega \left(\sup_{\substack{\bar{w} \in \bar{\chi}_w \\ z^\perp \in \chi^\perp}} (G^\perp)^\top G^\perp\right) dx \leq \beta^2 I_{n \times n}, \quad \forall t \geq 0. \quad (21)$$

4.
$$\lambda_1 \lambda_2 > \sigma^2, \quad \text{where} \quad \sigma := \frac{\beta}{2\sqrt{m_{1,*} m_{2,*}}}. \quad (22)$$

*Here,*

$$G := \left(\frac{\partial g_1}{\partial z^\perp}\right)^\top M_1 + M_2 \frac{\partial g_2}{\partial \bar{w}}, \quad G^\perp := G - \int_\Omega G \, dx, \quad (23)$$

and $m_{1,*}, m_{2,*} > 0$ are constants, guaranteed to exist by the uniform positive definiteness assumption, such that
$$m_{1,*}I_{n,n} \leq M_1(t), \quad \forall t \geq 0,$$
and
$$m_{2,*}I_{m,m} \leq M_2(t,x), \quad \forall t \geq 0, \ \forall x \in \Omega.$$

Then, the system is contracting with contraction rate
$$\lambda^* := -\left[\frac{\lambda_1 + \lambda_2}{2} - \sqrt{\left(\frac{\lambda_1 - \lambda_2}{2}\right)^2 + \sigma^2}\right].$$

*Proof.* Let $[\bar{w}_1, z_1^\perp]^\top \in \bar{\chi}_w \times \chi^\perp$ and $[\bar{w}_2, z_2^\perp]^\top \in \bar{\chi}_w \times \chi^\perp$ be any two solutions to (18), and define the error variables as $\bar{e} = \bar{w}_1 - \bar{w}_2$ and $e^\perp = z_1^\perp - z_2^\perp$. We define the error metric as
$$v(t) = \bar{e}^\top M_1 \bar{e} + \int_\Omega e^{\perp,\top} M_2 e^\perp \, dx,$$
where $v(t)$ represents a measure of the difference between the two solutions. Contraction is equivalent to $v(t)$ satisfying the inequality
$$v(t) \leq v(0)e^{-\lambda^* t}.$$
By taking the time derivative of $v(t)$ and applying the theorem's conditions, this inequality can be established. The full proof is provided in the SI Section 7.5. □

*Remark* 1. The condition given by (19) guarantees the contraction of the isolated subsystem described by
$$\frac{d\bar{w}(t)}{dt} = f_1(t, \bar{w}), \tag{24}$$
and (20) guarantees the contraction of the isolated subsystem described by
$$\frac{\partial z^\perp(t,x)}{\partial t} = \mathcal{L}_\Theta(D, z^\perp) + f_2^\perp(t, x, z^\perp), \quad \forall x \in \Omega, \quad \forall t > 0, \tag{25}$$
$$J(\Theta^i, D^i, z^{\perp,i}) \cdot \nu = 0, \quad x \in \partial\Omega, \quad t > 0, \quad i = 1, \ldots, m. \tag{26}$$

Conditions (21) and (22) handle the interconnection between these two contracting subsystems. This is analogous to the well-established feedback combination of two contracting subsystems [25].

*Remark* 2. The use of a virtual system allows us to extend the class of systems to which Theorem 1 can be applied to. By showing that a virtual system is contracting toward a particular solution, and that the original system is a particular solution of the virtual system, we can conclude that the original system will also approach the particular solution of the virtual system.

Next, we utilize a *virtual system* to expand the class of systems to which we can apply Theorem 1. Specifically, if we can show that a virtual system is contracting, that it approaches a particular solution, and the original system is a solution of the virtual system, then, by the definition of contraction, all trajectories of the virtual system will converge exponentially fast toward each other. Consequently, the original system will also approach the particular solution of the virtual system.

**Corollary** 1. *Consider the system*
$$\frac{d\bar{w}(t)}{dt} = f_1(t, \bar{w}) + \int_\Omega g_1(t, x', \bar{w}, z^\perp) \, dx', \quad \forall t > 0,$$
$$\frac{\partial z^\perp(t,x)}{\partial t} = \mathcal{L}_\Theta(D, z^\perp) + f_2^\perp(t, x, \bar{w}, z^\perp) + g_2^\perp(t, x, \bar{w}), \quad x \in \Omega, \quad t > 0, \tag{27}$$
$$J(\Theta^i, D^i, z^{\perp,i}) \cdot \nu = 0, \quad x \in \partial\Omega, \quad t > 0, \quad i = 1, \ldots, m,$$

and assume that $\bar{w}(t) = 0$ and $z^\perp(t, x) = 0$ are solutions to the system, that is,

$$
\begin{aligned}
\int_\Omega g_1(t, x, \bar{w}, 0)\, dx &= 0, \quad \forall t \geq 0,\ \forall \bar{w} \in \bar{\chi}_w, \\
f_2^\perp(t, x, \bar{w}, 0) &= 0, \quad \forall t \geq 0,\ x \in \Omega,\ \bar{w} \in \bar{\chi}_w, \\
f_1(t, 0) &= 0, \quad \forall t \geq 0, \\
g_2^\perp(t, x, 0) &= 0, \quad \forall t \geq 0,\ x \in \Omega,\ \hat{\bar{w}}(t) \in \bar{\chi}_w.
\end{aligned} \tag{28}
$$

Next, we define the following virtual system, where $\bar{w}(t)$ can be considered as a time-varying input:

$$
\begin{aligned}
\frac{d\bar{y}_w(t)}{dt} &= f_1(t, \bar{y}_w) + \int_\Omega g_1(t, x', \bar{w}, y_z^\perp)\, dx', \quad t > 0, \\
\frac{\partial y_z^\perp(t, x)}{\partial t} &= \mathcal{L}_\Theta(D, y_z^\perp) + f_2^\perp(t, x, \bar{w}, y_z^\perp) + g_2^\perp(t, x, \bar{y}_w), \quad x \in \Omega,\ t > 0, \\
J(\Theta^i, D^i, y_{z,i}^\perp) \cdot \nu &= 0, \quad x \in \partial\Omega,\ t > 0,\ i = 1, \ldots, m.
\end{aligned} \tag{29}
$$

We observe that $\bar{y}_w(t) = \bar{w}(t) = 0$ and $y_z^\perp(t, x) = z^\perp(t, x) = 0$ are solutions to the virtual system (29). Furthermore, the virtual system is in the form stipulated by equation (18), allowing us to apply Theorem 1 directly to demonstrate contraction. Furthermore, contraction implies that $\bar{w}(t) \to 0$ and $z^\perp(t, x) \to 0$.

**Corollary** 2. *In addition to the conditions of Corollary 1, suppose that **either** one of the following conditions holds for (29):*

1. $\int_\Omega g_1(t, x', \bar{w}, y_z^\perp)\, dx' \equiv 0$ and $\frac{\partial g_2}{\partial \bar{y}_w}$ is uniformly bounded,

2. $g_2^\perp(t, x, \bar{y}_w) \equiv 0$ and $\frac{\partial g_1}{\partial y_z^\perp}$ is uniformly bounded,

*for all $t \geq 0$, $x \in \Omega$, $\bar{w}, \bar{y}_w \in \bar{\chi}_w$, and $z^\perp, y_z^\perp \in \chi^\perp$. Then, to establish contraction of the system in (29), it suffices to verify conditions (19) and (20) from Theorem 1. Under these assumptions, $\bar{y}_w(t) \to 0$ and $y_z^\perp(t, x) \to 0$ with a contraction rate given by $\lambda^* = \min(\lambda_1, \lambda_2)$.*

*Proof.* If condition (1) is satisfied, then the dynamics of $\bar{y}_w$ are decoupled from $y_z^\perp$ in (29), and thus (19) implies contraction of the $\bar{y}_w$ dynamics. Thereafter, $\bar{y}_w$ can be treated as a time-varying input to the $y_z^\perp$ dynamics. Consequently, (20) and the uniform boundedness of $\frac{\partial g_2}{\partial \bar{y}_w}$ together imply contraction of the $y_z^\perp$ dynamics.

Conversely, if condition (2) is satisfied, then the $y_z^\perp$ dynamics are decoupled from $\bar{y}_w$ in (29), and thus (20) implies contraction of the $y_z^\perp$ dynamics. Thereafter, $y_z^\perp$ can be treated as a time-varying input to the $\bar{y}_w$ dynamics. Consequently, (19) and the uniform boundedness of $\frac{\partial g_1}{\partial y_z^\perp}$ together imply contraction of the $\bar{y}_w$ dynamics.

This framework is an example of a hierarchically combined system of two contracting dynamics in contraction theory [20, 26]. □

## 5 Examples

In this section, we revisit the examples from Section 3 and apply the results from Section 4 to analyze them. Additionally, we extend our analysis to more general cases, further demonstrating the applicability of our theoretical findings

### 5.1 Scalar State with Fickian Diffusion and Spatially Varying Reaction Terms

Consider a scalar state variable $z(t, x) \in \mathbb{R}$ defined on the one-dimensional spatial domain $\Omega = (0, 1)$ with dynamics described by

$$
\begin{aligned}
\frac{\partial z(t, x)}{\partial t} &= d\frac{\partial^2 z(t, x)}{\partial x^2} + a(x) z(t, x), \quad x \in \Omega,\ t > 0, \\
\frac{\partial z(t, x)}{\partial x} &= 0, \quad x \in \partial\Omega,\ t > 0,
\end{aligned} \tag{30}
$$

where $d > 0$ is the constant diffusion coefficient, and $a(x) : \Omega \to \mathbb{R}$ is a smooth function representing a spatially varying reaction rate.

Define the spatial average $\bar{w}(t) := \int_\Omega z(t, x) \, dx$ and the deviation from the average $z^\perp(t, x) := z(t, x) - \bar{w}(t)$. Taking the time derivative of both, we obtain

$$\frac{d\bar{w}(t)}{dt} = \bar{a}\bar{w}(t) + \int_\Omega a(x')z^\perp(t, x') \, dx', \tag{31}$$

$$\frac{\partial z^\perp(t, x)}{\partial t} = d\frac{\partial^2 z^\perp(t, x)}{\partial x^2} + \left(a(x)z^\perp(t, x) - \int_\Omega a(x')z^\perp(t, x') \, dx'\right) + a^\perp(x)\bar{w}(t), \quad x \in \Omega, \quad t > 0,$$

$$\frac{\partial z^\perp(t, x)}{\partial x} = 0, \quad x \in \partial\Omega, \quad t > 0,$$

where $\bar{a} := \int_\Omega a(x) \, dx$ is the spatial average of $a(x)$ and $a^\perp(x) = a(x) - \bar{a}$ it's the deviation from the spatial average.

This system fits into the framework of Theorem 1 with the following identifications:

$$\begin{aligned} f_1(t, \bar{w}) &= \bar{a}\bar{w}(t), & g_1(t, x, z^\perp) &= a(x)z^\perp(t, x), \\ f_2(t, x, z^\perp) &= a(x)z^\perp(t, x), & g_2(t, x, \bar{w}) &= a^\perp(x)\bar{w}(t). \end{aligned} \tag{32}$$

Note that $\Psi_{\Theta,D}(x) = 1$ since $\theta = 1/2$ and $d(x) = d$ is constant. Next, we apply Theorem 1 to (32) with $M_1 = M_2 = 1$ and defining $a^* := \sup_{x \in \Omega} a(x)$:

1. *Stability of the Spatial Average Dynamics* (19): We have that $\frac{\partial f_1}{\partial \bar{w}} = \bar{a}$. Thus, the inequalities,

$$\lambda_1 > 0 \quad \text{and} \quad \frac{1}{2}\left[\bar{a} + \bar{a}^T\right] \leq -\lambda_1, \tag{33}$$

   are satisfied provided that $\bar{a}$ is strictly negative, in which case we may take $\lambda_1 = -\bar{a}$.

2. *Stability of the Spatial Deviation Dynamics* (20): We have $\frac{\partial f_2}{\partial z^\perp} = a(x)$ and $\Lambda_{\Theta,D} = d\pi^2$. Thus, the inequalities

$$\lambda_2 > 0 \quad \text{and} \quad \frac{1}{2}\left[(a(x) - d\pi^2) + (a^T(x) - d\pi^2)\right] \leq -\lambda_2, \quad \forall x \in \Omega, \tag{34}$$

   are satisfied provided that diffusion is sufficiently fast (i.e., $d\pi^2 > a^*$), in which case we may take $\lambda_2 = d\pi^2 - a^*$.

3. *Boundedness of the Coupling Terms* (21): We have that $\frac{\partial g_1}{\partial z^\perp} = a(x)$ and $\frac{\partial g_2}{\partial \bar{w}} = a^\perp(x)$. Thus, the innequlity

$$\int_\Omega \left(a^\perp(x) + a^{\perp,\top}(x)\right)^2 dx \leq \beta^2, \tag{35}$$

   is satisfied for $\beta^2 = 4\|a^\perp(x)\|^2$.

4. *Small Gain Condition* (22):

$$\lambda_1 \lambda_2 > \frac{\beta^2}{4} \implies |\bar{a}|(d\pi^2 - a^*) \geq \|a^\perp(x)\|^2 \implies d\pi^2 > \frac{\|a^\perp(x)\|^2}{|\bar{a}|} + a^* \tag{36}$$

Altogether, the system (30) is contracting if the space-average dynamics are negative ($\bar{a} < 0$) and the diffusion rate sufficiently exceeds the spatial heterogeneity in $a(x)$, namely, $d\pi^2 > \frac{\|a^\perp(x)\|^2}{|\bar{a}|} + a^*$. Note that these conditions are independent of the state variables, so $\bar{\chi}_w$ and $\chi^\perp$ may be taken as the entire real line.

*Remark 3. The contraction of the linear system (30) is equivalent to the operator*

$$\mathcal{L}_a(z) := -d\frac{\partial^2 z}{\partial x^2} - a(x)z$$

*being positive definite under Neumann boundary conditions. Previous results (e.g., [27]) suggest that a sufficient condition for positive definiteness is $a(x) < 0$ for all $x \in \Omega$. However, our result provides a less restrictive condition that allows $a(x)$ to be positive on subsets of $\Omega$.*

## 5.2 Revisiting Example 3.1

We reconsider the dynamical system (1) and verify the inequalities in Section 5.1 with

$$d = \frac{\epsilon}{\pi^2}, \quad a(x) = -\epsilon + \sin(\omega x) - \int_0^1 \sin(\omega x') \, dx'.$$

The inequalities (33) hold under $\epsilon > 0$, giving

$$\lambda_1 = -\bar{a} = \epsilon.$$

To ensure $\lambda_2 > 0$ in (34), we require

$$\lambda_2 = \epsilon - \sup_{x \in (0,1)} \left[ -\epsilon + \sin(\omega x) - \int_0^1 \sin(\omega x') \, dx' \right] > 0.$$

For simplicity, assume $\omega \ll 1$, so $\sin(\omega x) \approx \omega x$ for $x \in (0, 1)$. Then

$$\lambda_2 \approx 2\epsilon - \frac{\omega}{2},$$

and positivity requires $\omega < 4\epsilon$.

Next, observe that

$$a^\perp(x) = \sin(\omega x) - \frac{1 - \cos(\omega)}{\omega},$$

so (35) is satisfied with

$$\beta = 4 \int_0^1 \left( \sin(\omega x) - \frac{1 - \cos(\omega)}{\omega} \right)^2 dx,$$

and for $\omega \ll 1$, $\beta \approx \frac{\omega^2}{3}$.

Finally, (36) is equivalent to

$$\epsilon \left( 2\epsilon - \frac{\omega}{2} \right) > \frac{\omega^2}{12},$$

giving

$$\omega < (\sqrt{33} - 3) \epsilon.$$

Hence, (1) is contracting if $\epsilon > 0$ and $\omega < (\sqrt{33} - 3)\epsilon$. This aligns with numerical simulations (Fig. 1-b), where for $\epsilon = 1 \times 10^{-2}$, the bound is about $\omega \lesssim 3 \times 10^{-2}$.

## 5.3 High-Dimensional Linear Systems with General Diffusion and Spatially Varying Reactions

We extend the system from equation (30) to higher dimensions, where $z(t, x) \in \mathbb{R}^n$, and incorporate time- and space-varying reaction terms using the general diffusion model defined earlier. For $\Theta \in \mathbb{R}^n$ and $D(x) \in \mathbb{R}^n$ such that $\Theta^i \in [0, 1]$, $D^i(x) > 0$, and $D^i(x)$ are smooth for all $i = 1, \ldots, n$, the system is described by:

$$\frac{\partial z(t, x)}{\partial t} = \mathcal{L}_\Theta(D, z) + A(t, x)z(t, x), \tag{37}$$
$$J(\Theta^i, D^i, z^i) \cdot \nu = 0, \quad \forall x \in \partial\Omega, \quad \forall t > 0, \quad \forall i = 1, \ldots, n,$$

where, $A(t, x) \in \mathbb{R}^{n \times n}$ is a smoothly varying reaction matrix, $\mathcal{L}_\Theta^i(D, z^\perp) = L(\Theta^i, D^i(x), z^{\perp, i})$, as given by equation (13), and $J(\Theta^i, D^i, z^i)$ as given by equation (12).

Let $\bar{w}(t) := \int_\Omega z(t, x) \, dx$ denote the spatial average, and define the spatially weighted deviations from the mean as:

$$z^\perp(t, x) := z(t, x) - \Psi_{\Theta, D}(x)\bar{w}(t),$$

where $\Psi_{\Theta, D}(x) = \mathrm{diag}\left(\psi_{\Theta^1, D^1}(x), \ldots, \psi_{\Theta^n, D^n}(x)\right)$, with each $\psi_{\Theta^i, D^i}(x)$ defined in equation (15).

The dynamics of $\bar{w}(t)$ and $z^\perp(t,x)$ are then:

$$\frac{d\bar{w}(t)}{dt} = \bar{A}_\Psi(t)\bar{w}(t) + \int_\Omega A(t,x) z^\perp(t,x)\, dx, \tag{38}$$

$$\frac{\partial z^\perp(t,x)}{\partial t} = \mathcal{L}_\Theta(D, z^\perp) + A(t,x)\left[z^\perp(t,x) + \Psi_{\Theta,D}(x)\bar{w}(t)\right]$$

$$- \Psi_{\Theta,D}(x) \int_\Omega A(t,x)\left[z^\perp(t,x) + \Psi_{\Theta,D}(x)\bar{w}(t)\right]\, dx,$$

$$J(\Theta^i, D^i, z^{\perp,i}) \cdot \nu = 0, \quad \forall x \in \partial\Omega, \quad \forall t > 0, \quad \forall i = 1, \ldots, n,$$

where

$$\bar{A}_\Psi(t) := \int_\Omega A(t,x) \Psi_{\Theta,D}(x)\, dx.$$

This system aligns with the structure required by Theorem 1, with the following correspondences:

$$f_1 = \bar{A}_\Psi(t)\bar{w}(t), \quad g_1 = A(t,x)z^\perp(t,x), \quad f_2 = A(t,x)z^\perp(t,x), \quad g_2 = A(t,x)\Psi_{\Theta,D}(x)\bar{w}(t). \tag{39}$$

To ensure contraction, we require the existence of positive definite matrices $M_1(t)$, $\Gamma(t)$, and constants $\lambda_1, \lambda_2, \beta > 0$ such that:

1. *Stability of the Spatial Average Dynamics (19):*

$$\frac{1}{2}\left[\dot{M}_1(t) + M_1(t)\bar{A}_\Psi(t) + \bar{A}_\Psi^T(t) M_1(t)\right] \leq -\lambda_1 M_1(t), \quad \forall t \geq 0. \tag{40}$$

This condition is equivalent to the ODE contraction condition [20] for the isolated $\bar{w}(t)$ dynamics in (38) when $z^\perp = 0$.

2. *Stability of the Spatial Deviation Dynamics (20):*

$$\frac{1}{2}\left[\dot{M}_2(t,x) + M_2(t,x)\Big(A(t,x) - \Lambda_{\Theta,D}\Big) + \Big(A^T(t,x) - \Lambda_{\Theta,D}\Big)M_2(t,x)\right] \leq -\lambda_2 M_2(t,x), \quad \forall t \geq 0, \forall x \in \Omega, \tag{41}$$

where $\Lambda_{\Theta,D}$ is the diagonal matrix of eigenvalues associated with the diffusion operator $\mathcal{L}_\Theta(D, z)$ (as defined in (16)). Here, $M_2(t,x) = \Gamma(t)\Psi_{\Theta,D}^{-1}(x)$ with $\Gamma(t)$ a positive definite diagonal matrix. This condition is equivalent to the contraction condition for the isolated $z^\perp(t,x)$ dynamics in (38) when $\bar{w}(t) = 0$.

3. *Boundedness of the Coupling Terms (21):*

$$\int_\Omega G^{T,\perp}(t,x) G^\perp(t,x)\, dx \leq \beta^2 I, \quad \forall t \geq 0,$$

where

$$G(t,x) = A^T(t,x)M_1(t) + \Gamma(t)\Psi_{\Theta,D}^{-1}(x)A(t,x)\Psi_{\Theta,D}(x), \quad G^\perp(t,x) = G(t,x) - \int_\Omega G(t,x)\, dx.$$

This condition controls the coupling between the dynamics of $\bar{w}(t)$ and $z^\perp(t,x)$.

In addition to the three condition above, to ensure contraction, the small gain conditions (22) also needs to be satisfied. The existence of $M_1(t)$ and $\lambda_1$ imposes a stability condition on the weighted spatial average matrix $\bar{A}_\Psi(t)$. Although the unweighted average $\bar{A}(t) := \int_\Omega A(t,x)\, dx$ might be stable or unstable, the weighted average with $\Psi_{\Theta,D}(x)$ can influence the overall stability of the dynamics.

The matrix $\Gamma(t)$ is chosen to ensure the diagonal stability of the operator $A(t,x) - \Lambda_{\Theta,D}$. The coupling term $G^\perp(t,x)$ captures the influence of spatial heterogeneity in $A(t,x)$ and the heterogeneity in the similarity transformation $\Psi_{\Theta,D}^{-1}(x)A(t,x)\Psi_{\Theta,D}(x)$ due to both diffusion and reaction dynamics.

By satisfying the above conditions, we can apply Theorem 1 to conclude that the system (37) is contracting. Consequently, this implies $z^\perp(t,x) \to 0$ and $\bar{w}(t) \to 0$, ensuring the system's stability. Notice that the conditions are independent of the state variables and thus $\bar{\chi}_w$ and $\chi^\perp$ can both be taken to be $\mathbb{R}^n$

## 5.4 Revisiting Example 3.2

We revisit the example described by the system (6). From (39), the system terms are given by:

$$f_1 = A\bar{z}(t), \quad g_1 = Az^\perp(t,x), \quad f_2 = Az^\perp(t,x), \quad g_2 = A\Psi_{\Theta,D}(x)\bar{z}(t),$$

where $\Theta = \begin{bmatrix} 1/2 & 1 \end{bmatrix}^T$, $D = \zeta \begin{bmatrix} 10 & 1/4v_r(x) \end{bmatrix}^T$, and $A = \begin{bmatrix} -1 & 1 \\ -1 & \frac{1}{2} \end{bmatrix}$.

We observe that:

$$\int_\Omega g_1\, dx = A \int_\Omega z^\perp(t,x)\, dx = 0.$$

Thus, by Corollary 2, to establish contraction of (6), it suffices to satisfy the first two conditions, (40)–(41).

The existence of $M_1$ and $\lambda_1$ in (40) is guaranteed by the fact that $A$ is Hurwitz, meaning the eigenvalues of $A$ have negative real parts $-0.25 \pm 0.66j$. The eigenvalue matrix of the differential operator, $\Lambda_{\Theta,D}$, is given by:

$$\Lambda_{\Theta,D} = \frac{\zeta}{\pi^2} \begin{bmatrix} 10 & 0 \\ 0 & \frac{1}{4} \end{bmatrix} \begin{bmatrix} \pi^2 & 0 \\ 0 & \pi^2 \nu \end{bmatrix}, \quad \text{where} \quad \nu = \frac{\min_{x \in (0,1)} v_r^2(x)}{\max_{x \in (0,1)} v_r(x)}. \tag{42}$$

The bound for $\nu$ was derived using (17). To satisfy the second condition of the theorem (41), we need to certify the diagonal stability of

$$\Psi_{\Theta,D}^{-1}(x)\left[A - \Lambda_{\Theta,D}\right] = \begin{bmatrix} 1 & 0 \\ 0 & v_r^{-1}(x) \end{bmatrix} \left[ \begin{bmatrix} -1 & 1 \\ -1 & \frac{1}{2} \end{bmatrix} - \zeta \begin{bmatrix} 10 & 0 \\ 0 & \frac{1}{4} \end{bmatrix} \begin{bmatrix} 1 & 0 \\ 0 & \nu \end{bmatrix} \right].$$

For a general $2 \times 2$ matrix $A$, diagonal stability is equivalent to the negative of matrix having all positive principle minors [28]. Therefore, for contraction of (2), we need that

$$(1 + 10\zeta) > 0, \ (-1/2 + 1/4\nu\zeta) > 0, \ \frac{(1+10\zeta)(-1/2 + 1/4\nu\zeta)}{\hat{v}_r(x)} + \frac{1}{\hat{v}_r(x)} > 0,$$

altogether this implies that

$$\zeta > \frac{2}{\nu}. \tag{43}$$

By (42), note that $\nu$ captures spatial heterogeneity. The approximated upper bound, shown as the blue curve in Fig. 2(c), indeed serves as an upper bound to the numerical estimate of $\zeta_{cr}$.

## 5.5 Revisiting Example 3.3

We revisit the biochemical system described by Equation (7), which models cellular translational dynamics, and demonstrate how to use the tools from this work to justify the approximations in Equations (9)–(10). Normalizing time by $\tau = t(\kappa + d)$, the spatiotemporal concentration dynamics become

$$\begin{aligned}
\frac{\partial m}{\partial \tau} &= L\left(1, \chi_m v_{r_m}, m\right) - \left(\frac{mR}{K} - c_r\right), \\
\frac{\partial R}{\partial \tau} &= L\left(1, \chi_r v_{r_r}, R\right) - \left(\frac{mR}{K} - c_r\right), \\
\frac{\partial c_r}{\partial \tau} &= L\left(1, \chi_c v_{r_c}, c_r\right) + \left(\frac{mR}{K} - c_r\right),
\end{aligned} \tag{44}$$

where $\chi_m = \frac{D_m}{\kappa+d}$, $\chi_r = \frac{d}{\kappa+d}$, $\chi_c = \frac{D_c}{\kappa+d}$, and $L$ is the diffusion operator defined in Equation (13), along with the no-flux boundary conditions Equation (14).

Integrating both sides of (44) over $\Omega$ and applying the boundary conditions, the space-averaged dynamics are described by

$$\frac{d\bar{m}}{d\tau} = -\int_\Omega \left(\frac{mR}{K}\right) dx + \bar{c}_r, \quad \frac{d\bar{R}}{d\tau} = -\int_\Omega \left(\frac{mR}{K}\right) dx + \bar{c}_r, \quad \frac{d\bar{c}_r}{d\tau} = \int_\Omega \left(\frac{mR}{K}\right) dx - \bar{c}_r. \tag{45}$$

]

From (45), we observe that the spaced-averaged concentration of total mRNA and total ribosomes in the cell are conserved quantities:

$$\frac{d}{d\tau}(\bar{c}_r + \bar{m}) = 0 \implies \bar{c}_r(\tau) + \bar{m}(\tau) = \bar{m}_T := \bar{c}_r(0) + \bar{m}(0) = \text{constant}, \quad (46)$$

$$\frac{d}{d\tau}(\bar{c}_r + \bar{R}) = 0 \implies \bar{c}_r(\tau) + \bar{R}(\tau) = \bar{R}_T := \bar{c}_r(0) + \bar{R}(0) = \text{constant}.$$

Consequently, to fully specify the space-averaged dynamics, it is sufficient to consider $\bar{c}_r(\tau)$. Next, we define the following error variables:

$$\bar{e}(t) = \bar{c}_r(t) - \theta \frac{\bar{m}(t)\bar{R}(t)}{K}, \quad m^\perp(\tau, x) = m(\tau, x) - \hat{v}_{r_m}(x)\bar{m}(\tau),$$
$$R^\perp(\tau, x) = R(\tau, x) - \hat{v}_{r_r}(x)\bar{R}(\tau), \quad c_r^\perp(\tau, x) = c_r(\tau, x) - \hat{v}_{r_c}(x)\bar{c}_r(\tau). \quad (47)$$

When the error variables in (47) are zero, the approximations in Equations (9)–(10) are satisfied, as desired. The dynamics of the error variables are given by

$$\frac{d\bar{e}}{d\tau} = \left(1 + \theta \frac{(\bar{m} + \bar{R})}{K}\right)\left[-\bar{e} + \frac{1}{K}\int_\Omega \left(\hat{v}_{r_m}(x)\bar{m}R^\perp + m^\perp \hat{v}_{r_r}(x)\bar{R} + m^\perp R^\perp\right) dx\right], \quad (48)$$

$$\frac{\partial m^\perp}{\partial \tau} = L\left(1, \chi_m v_{r_m}, m^\perp\right) - \frac{1}{K}\left(\hat{v}_{r_m}(x)\bar{m}R^\perp + m^\perp \hat{v}_{r_r}(x)\bar{R} + m^\perp R^\perp\right) + c_r^\perp + \hat{v}_{r_c}(x)\bar{e} + \hat{v}_{r_m}(x)\frac{d\bar{e}/d\tau}{1 + \theta\frac{\bar{m}+\bar{R}}{K}},$$

$$\frac{\partial R^\perp}{\partial \tau} = L\left(1, \chi_r v_{r_r}, R^\perp\right) - \frac{1}{K}\left(\hat{v}_{r_m}(x)\bar{m}R^\perp + m^\perp \hat{v}_{r_r}(x)\bar{R} + m^\perp R^\perp\right) + c_r^\perp + \hat{v}_{r_c}(x)\bar{e} + \hat{v}_{r_r}(x)\frac{d\bar{e}/d\tau}{1 + \theta\frac{\bar{m}+\bar{R}}{K}},$$

$$\frac{\partial c_r^\perp}{\partial \tau} = L\left(1, \chi_c v_{r_c}, c_r^\perp\right) + \frac{1}{K}\left(\hat{v}_{r_m}(x)\bar{m}R^\perp + m^\perp \hat{v}_{r_r}(x)\bar{R} + m^\perp R^\perp\right) - c_r^\perp - \hat{v}_{r_c}(x)\bar{e} - \hat{v}_{r_c}(x)\frac{d\bar{e}/d\tau}{1 + \theta\frac{\bar{m}+\bar{R}}{K}}.$$

To derive (48), we used $v_{r_c}(x) = v_r(x)v_{r_m}(x)$, which follows from the factorization of the polysome available volume into independent contributions from ribosome-free mRNA and ribosomes, as shown in Equation (S12) of [2]. To transform (48) into a form to which Theorem 1 is applicable, we define the following virtual dynamics:

$$\frac{d\bar{y}_e}{d\tau} = \left(1 + \theta \frac{(\bar{m} + \bar{R})}{K}\right)\left[-\bar{y}_e + \frac{1}{K}\int_\Omega \left(\hat{v}_{r_m}(x)\bar{m}y_R^\perp + y_m^\perp \hat{v}_{r_r}(x)\bar{R} + \frac{1}{2}\left(R^\perp y_m^\perp + y_R^\perp m^\perp\right)\right) dx\right], \quad (49)$$

$$\frac{\partial y_m^\perp}{\partial \tau} = L\left(1, \chi_m v_{r_m}, y_m^\perp\right) - \frac{1}{K}\left(\hat{v}_{r_m}(x)\bar{m}y_R^\perp + y_m^\perp \hat{v}_{r_r}(x)\bar{R} + \frac{1}{2}\left(R^\perp y_m^\perp + y_R^\perp m^\perp\right)\right) + y_{c_r}^\perp + \hat{v}_{r_c}(x)\bar{y}_e + \hat{v}_{r_m}(x)\frac{d\bar{y}_e/d\tau}{1 + \theta\frac{\bar{m}+\bar{R}}{K}},$$

$$\frac{\partial y_R^\perp}{\partial \tau} = L\left(1, \chi_r v_{r_r}, y_R^\perp\right) - \frac{1}{K}\left(\hat{v}_{r_m}(x)\bar{m}y_R^\perp + y_m^\perp \hat{v}_{r_r}(x)\bar{R} + \frac{1}{2}\left(R^\perp y_m^\perp + y_R^\perp m^\perp\right)\right) + y_{c_r}^\perp + \hat{v}_{r_c}(x)\bar{y}_e + \hat{v}_{r_r}(x)\frac{d\bar{y}_e/d\tau}{1 + \theta\frac{\bar{m}+\bar{R}}{K}},$$

$$\frac{\partial y_{c_r}^\perp}{\partial \tau} = L\left(1, \chi_c v_{r_c}, y_{c_r}^\perp\right) + \frac{1}{K}\left(\hat{v}_{r_m}(x)\bar{m}y_R^\perp + y_m^\perp \hat{v}_{r_r}(x)\bar{R} + \frac{1}{2}\left(R^\perp y_m^\perp + y_R^\perp m^\perp\right)\right) - y_{c_r}^\perp - \hat{v}_{r_c}(x)\bar{y}_e - \hat{v}_{r_c}(x)\frac{d\bar{y}_e/d\tau}{1 + \theta\frac{\bar{m}+\bar{R}}{K}},$$

where $\bar{m}(\tau), \bar{R}(\tau), R^\perp(\tau, x)$, and $m^\perp(\tau, x)$ are interpreted as time-varying inputs to the system, and consequently the virtual dynamics are linear. Notice that $\bar{y}_e = y_m^\perp = y_R^\perp = y_{c_r}^\perp = 0$ and $\bar{y}_e = \bar{e}$, $y_m^\perp = m^\perp$, $y_R^\perp = R^\perp$, $y_{c_r}^\perp = c_r^\perp$ are solutions to (49).

The virtual system (49) has the form stipulated by Equation (29) in Corollary 1, with $\bar{y}_w = \bar{y}_e$, $y_z^\perp = \left[y_m^\perp, \; y_R^\perp, \; y_{c_r}^\perp\right]^T$, and:

$$f_1(\tau, \bar{y}_w) = -\eta(\tau)\bar{y}_w, \quad g_1(\tau, x, y_z^\perp) = \eta(\tau)u^T(\tau, x)y_z^\perp(\tau, x), \quad g_2(\tau, x, \bar{y}_w) = -\hat{v}_{r_c}(x)v\bar{y}_w,$$
$$f_2(\tau, x, y_z^\perp) = v\tilde{u}^T(\tau, x)y_z^\perp(\tau, x), \quad D = \begin{bmatrix}\chi_m v_{r_m}, & \chi_r v_{r_r}, & \chi_c v_{r_c}\end{bmatrix}^T, \quad \Theta = \begin{bmatrix}1, & 1, & 1\end{bmatrix}^T,$$

where:

$$\eta(\tau) = 1 + \theta\frac{(\bar{m}(\tau) + \bar{R}(\tau))}{K},$$

$$u^T(\tau,x) = \left[\frac{\hat{v}_{r_r}(x)\bar{R}(\tau) + \frac{1}{2}R^\perp(\tau,x)}{K}, \quad \frac{\hat{v}_{r_m}(x)\bar{m}(\tau) + \frac{1}{2}m^\perp(\tau,x)}{K}, \quad 0\right],$$

$$v^T = \begin{bmatrix} -1, & -1, & 1 \end{bmatrix},$$

$$\tilde{u}^T(\tau,x) = u^T(\tau,x) + \begin{bmatrix} 0, & 0, & -1 \end{bmatrix}.$$

Next, we verify the contraction conditions for (49). Given that the virtual dynamics are linear, we take $\bar{\chi}_w = \mathbb{R}$ and $\chi^\perp = \mathbb{R}^3$.

1. *Stability of the Spatial Average Dynamics* (19): We have $\frac{\partial f_1}{\partial \bar{y}_w} = -\eta(\tau)$. By letting $M_1 = 1$, we obtain an approximate bound for $\lambda_1$ from the inequality

$$\frac{1}{2}\left[-\eta(\tau) - \eta^T(\tau)\right] \leq -\lambda_1, \quad \forall \tau \geq 0,$$

which is satisfied with $\lambda_1 = 1$, leveraging the positivity of $\bar{m}(\tau)$ and $\bar{R}(\tau)$ (see SI Section 7.6).

2. *Stability of the Spatial Deviation Dynamics* (20): We have $\frac{\partial f_2}{\partial y_z^\perp} = v\tilde{u}^T(\tau,x)$ and letting $\Gamma(t) = I_{3,3}$, we approximate a bound for $\lambda_2$ in the inequality

$$\frac{1}{2}\left[\Psi_{\Theta,D}^{-1}(x)\left(v\tilde{u}^T(\tau,x) - \Lambda_{\Theta,D}\right) + \left(\tilde{u}(\tau,x)v^T - \Lambda_{\Theta,D}\right)\Psi_{\Theta,D}^{-1}(x)\right] \leq -\lambda_2 \Psi_{\Theta,D}^{-1}(x), \quad \forall \tau \geq 0, \forall x \in \Omega.$$

The inequality is satisfied for

$$\lambda_2 = \Lambda_{\Theta,D,*} - \sqrt{3\frac{\Psi_{\Theta,D}^*}{\Psi_{\Theta,D,*}}}\beta_u,$$

where

$$\Lambda_{\Theta,D,*} := \min_{i=\{1,2,3\}} \Lambda_{\Theta,D}^{i,i}, \quad \Psi_{\Theta,D,*} := \min_{i=\{1,2,3\}} (\inf_{x \in \Omega} \Psi_{\Theta,D}^{i,i}(x)),$$

$$\Psi_{\Theta,D}^* := \max_{i=\{1,2,3\}} (\sup_{x \in \Omega} \Psi_{\Theta,D}^{i,i}(x)), \quad \beta_u^2 := \left(1/2R^{\perp,*} + \hat{v}_{r_r}^* \bar{R}_T\right)^2 + \left(1/2m^{\perp,*} + \hat{v}_{r_m}^* \bar{m}_T\right)^2 + 1, \quad (50)$$

where $R^{\perp,*}$ and $m^{\perp,*}$ denote bounds on $R^\perp(\tau,x)$ and $m^\perp(\tau,x)$, respectively, valid for all $t \geq 0$ and $x \in \Omega$ and $\hat{v}_{r_m}^* := \sup_{x\in\Omega} \hat{v}_{r_m}(x)$ and $\hat{v}_{r_r}^* := \sup_{x\in\Omega} \hat{v}_{r_r}(x)$. The contraction rate $\lambda_2$ is strictly positive for sufficiently large $\Lambda_{\Theta,D,*}$, that is, for sufficiently fast diffusion rates. For full details see SI Section 7.6.

3. *Boundedness of the Coupling Terms* (21): We have $\frac{\partial g_1}{\partial y_z^\perp} = \eta(t)u^T(\tau,x)$ and $\frac{\partial g_2}{\partial \bar{y}_w} = -\hat{v}_{r_c}(x)v$ and thus from the coupling term equation (23) we have that

$$G(\tau,x) = \eta(\tau)u(\tau,x) + \psi_{\Theta,D}^{-1}(x)\hat{v}_{r_c}(x)v = \eta(\tau)u(\tau,x) + \begin{bmatrix} \frac{\hat{v}_{r_r}(x)}{\theta} & \frac{\hat{v}_{r_m}(x)}{\theta} & -1 \end{bmatrix}^T$$

$$G^\perp(\tau,x) = \eta(\tau)u'(\tau,x) + \begin{bmatrix} \frac{\hat{v}_{r_r}(x)-1}{\theta} & \frac{\hat{v}_{r_m}(x)-1}{\theta} & 0 \end{bmatrix}^T$$

where $u'(\tau,x) = \left[\frac{[\hat{v}_{r_r}(x)-1]\bar{R}(\tau) + \frac{1}{2}R^\perp(\tau,x)}{K}, \quad \frac{[\hat{v}_{r_m}(x)-1]\bar{m}(\tau) + \frac{1}{2}m^\perp(\tau,x)}{K}, \quad 0\right]^T$. Therefore,

$$\beta = \beta_{u,\eta} + \beta_h,$$

where

$$\beta_{u,\eta}^2 := \left(1 + \theta\frac{(\bar{m}_T + \bar{R}_T)}{K}\right)\int_\Omega \left(\left(\frac{[\hat{v}_{r_r}(x) - 1]\bar{R}_T + \frac{1}{2}R^{\perp,*}}{K}\right)^2 + \left(\frac{[\hat{v}_{r_m}(x) - 1]\bar{m}_T + \frac{1}{2}m^{\perp,*}}{K}\right)^2\right)dx$$

,

$$\beta_h^2 := \frac{1}{\theta^2}\int_\Omega \left((\hat{v}_{r_r}(x) - 1)^2 + (\hat{v}_{r_m}(x) - 1)^2\right)dx.$$

satisfies this inequality (SI Section 7.6 for full details).

4. *Small Gain Condition* (22):

$$\lambda_1 \lambda_2 > \frac{\beta^2}{4} \implies \left(\Lambda_{\Theta,D,*}\right) > \left(\frac{\beta_{u,\eta} + \beta_h}{2}\right)^2 + \sqrt{3 \frac{\Psi^*_{\Theta,D}}{\Psi_{\Theta,D,*}}} \beta_u,$$

Altogether, this implies that, for sufficiently fast diffusion, the dynamics (49) are contracting and the approximations (9)–(10) hold. In [6], these approximations were demonstrated assuming that diffusion was arbitrarily faster than the reaction dynamics—an assumption that is not always valid, particularly for the timescale associated with the binding/unbinding of ribosomes to mRNA [29]. Consequently, the results of this example provide a finite bound on the diffusion rate based on spatial heterogeneity in the reaction dynamics and diffusion.

Although, in the virtual system (49), we set $\bar{\chi}_w = \mathbb{R}$ and $\chi^\perp = \mathbb{R}^3$, the contraction conditions depend on bounds for the pointwise and space-averaged concentrations $m^{\perp,*}$, $R^{\perp,*}$, $\bar{m}_T$, and $\bar{R}_T$. From (46) and as shown in SI Section 7.6, these bounds are determined by the system's initial conditions. Consequently, the necessary speed of diffusion required to ensure contraction depends on the concentration range of interest.

*Remark 4. The structure of (48), where only a one-dimensional subset of the space-averaged dynamics is present, motivates the form of Equation (18). Instead of considering the full state space of $\bar{z}(t)$, we focus on $\bar{w}(t)$, which represents a subset of the space-averaged dynamics. This highlights an example of a partially contracting system, where convergence occurs within a subset of the state space while the overall system exhibits reduced-dimensional behavior [30].*

## 5.6 High-Dimensional Nonlinear Systems with General Diffusion and Spatially Varying Reactions

Next, we show how to transform a general nonlinear RD system into a form that Theorem 1 can be applied. Specifically, consider the a system identical to (37) but replace $A(t,x)z(t,x)$ for a general nonlinear function $f(t,x,z) \in \mathbb{R}^n$, the dynamics are described by,

$$\frac{\partial z(t,x)}{\partial t} = \mathcal{L}_\Theta(D,z) + f(t,x,z), \quad \forall x \in \Omega, \quad \forall t > 0, \tag{51}$$
$$J(\Theta^i, D^i, z^i) \cdot \nu = 0, \quad \forall x \in \partial\Omega, \quad \forall t > 0, \quad \forall i = 1,\ldots,n.$$

Let $\bar{z}(t) := \int_\Omega z(t,x)\,dx$ and $z^\perp(t,x) = z(t,x) - \Psi_{\Theta,D}(x)\bar{z}(t)$, the system can be equivalently written as

$$\frac{d\bar{z}(t)}{dt} = \int_\Omega f\big(t,x,\Psi_{\Theta,D}(x)\bar{z}(t) + z^\perp(t,x)\big)\,dx, \quad \forall t > 0,$$
$$\frac{\partial z^\perp(t,x)}{\partial t} = \mathcal{L}_\Theta(D,z^\perp) + f^\perp\big(t,x,\Psi_{\Theta,D}(x)\bar{z}(t) + z^\perp(t,x)\big), \quad \forall x \in \Omega, \quad \forall t > 0, \tag{52}$$
$$J(\Theta^i, D^i, z^{\perp,i}) \cdot \nu = 0, \quad \forall x \in \partial\Omega, \quad \forall t > 0, \quad \forall i = 1,\ldots,n.$$

Here,
$$f^\perp\big(t,x,\Psi_{\Theta,D}(x)\bar{z}(t) + z^\perp(t,x)\big) = f\big(t,x,\Psi_{\Theta,D}(x)\bar{z}(t) + z^\perp(t,x)\big)$$
$$- \Psi_{\Theta,D}(x) \int_\Omega f\big(t,x,\Psi_{\Theta,D}(x)\bar{z}(t) + z^\perp(t,x)\big)\,dx. \tag{53}$$

As in Corollary 1, we can define the following virtual system:

$$\frac{d\bar{y}_w(t)}{dt} = \underbrace{\int_\Omega f\big(t,x,\Psi_{\Theta,D}(x)\bar{y}_w(t)\big)\,dx}_{f_1(t,\bar{y}_w)}$$
$$+ \left[\int_\Omega \underbrace{f\big(t,x,\Psi_{\Theta,D}(x)\bar{z}(t) + y_z^\perp(t,x)\big) - f\big(t,x,\Psi_{\Theta,D}(x)\bar{z}(t)\big)}_{g_1(t,x,y_z^\perp)}\,dx\right], \quad \forall t>0,$$
$$\frac{\partial y_z^\perp(t,x)}{\partial t} = \mathcal{L}_\Theta(D,y_z^\perp) + \left[\underbrace{f^\perp\big(t,x,\Psi_{\Theta,D}(x)\bar{z}(t) + y_z^\perp(t,x)\big) - f^\perp\big(t,x,\Psi_{\Theta,D}(x)\bar{z}(t)\big)}_{f_2^\perp(t,x,y_z^\perp)}\right]$$
$$+ \underbrace{f^\perp\big(t,x,\Psi_{\Theta,D}(x)\bar{y}_w(t)\big)}_{g_2^\perp(t,x,\bar{y}_z)}, \quad \forall x\in\Omega,\quad \forall t>0, \qquad (54)$$
$$J(\Theta^i, D^i, y_{z,i}^\perp)\cdot \nu = 0, \quad \forall x\in\partial\Omega, \quad \forall t>0, \quad \forall i=1,\ldots,n.$$

Here, $\bar{z}(t)$ can be considered a time-varying input to (54). If the following relationship holds:

$$f(t,x,0) = 0, \quad \forall x\in\Omega, \quad \forall t>0,$$

then it implies that $\bar{y}_w = \bar{z} = 0$ and $y_z^\perp = z^\perp = 0$ are all solutions to (54), and the assumptions of Corollary 1 are satisfied. Thus, we can apply Theorem 1 to show contraction. The contraction conditions look identical to those in the linear space-varying example (Section 5.3) except that $A$ is replaced with $\frac{\partial f}{\partial z}$.

*Remark 5. The system described by (54) automatically satisfies the second condition stipulated by Corollary 2 if $\Psi_{\Theta,D}(x)$ and $f(t,x,z)$; are explicitly independent of $x$. Therefore, to check that $z^\perp \to 0$ with rate $\lambda_2$, it is only necessary to check (20) in Theorem 1 (see Corollary 2). This is in agreement with the results from [10, 11].*

## 6 Concluding remarks

Spatial heterogeneity can destabilize RD systems (Figures 1(a) and 2(b)) and can also lead to asymptotic solutions that approach heterogeneous spatial profiles (9). In this work, we employed a contraction approach to derive conditions that guarantee solutions to RD systems converge to each other exponentially fast, independent of initial conditions. The dynamics of the RD system are decomposed into the null space ($\Psi_{\Theta,D}(x)\bar{z}(t)$ in (52)) of the diffusion operator (13) and the residual dynamics ($z^\perp(t,x)$ in (52)). The first condition in Theorem 1 is identical to the contraction conditions of the isolated ODE system described by the weighted space-averaged dynamics in (24), which is similar to the approach used in [20]. The second condition in Theorem 1 is analogous to that in [10] and corresponds to the requirement that the isolated residual dynamics (25) be contracting. The key distinction between (20) and the condition presented in [10] is that the contraction metric $M_2$ has a spatial dependence via $\Psi_{\Theta,D}(x)$, ensuring the self-adjointness of the linear differential operator (13). Consequently, the condition must hold for all $x\in\Omega$, making it more stringent. Finally, conditions three and four in Theorem 1 handle the interconnection between (24) and (25) and are analogous to the results of feedback interconnected contracting system [25], which is a special case of the small gain theorem [31].

By using contraction theory, we can leverage tools such as partial contraction (Remark 4), virtual dynamics (Corollary 1), hierarchical combinations (Corollary 2), and feedback combinations (Remark 1) of contracting subsystems to simplify the analysis. The dynamical system (18) corresponding to Theorem 1 was presented with the space-averaged state variable $\bar{w}(t)$ instead of $\bar{z}(t)$ to account for partial contraction of the space-averaged dynamics, as observed in biochemical systems due to conservation laws (Remark 4). We applied the results of Section 4 to demonstrate contraction of both linear and nonlinear systems and revisited the examples from Section 3. In particular, we demonstrated the convergence of quasi-steady-state (QSS) solutions in biomolecular reactions occurring in the spatially heterogeneous cellular environment and how spatial heterogeneity modulates the binding rates in the space-averaged solutions (Section 5.5). This work extends the results of [6] by providing a finite bound on the diffusion rate required for the contraction of the QSS solution, relaxing the assumption of arbitrarily fast diffusion. This aligns with observed timescales in biology, where diffusion and biomolecular binding/unbinding can occur on similar timescales [29].

While the second condition in Theorem 1 requires the existence of a diagonal matrix $\Gamma(t)$, if the diffusion operators are identical for a subset of the states of $z^\perp(t,x)$, we can relax the diagonal requirement to only require a block-diagonal matrix, similar to the approach used in [11]. Furthermore, [2] demonstrated that the spatial distribution of the chromosomal mesh inside the cell varies with time as the cell divides. This implies that spatial heterogeneity varies over time, as modeled in [6]. Therefore, it is of interest to extend the results from this work to time-varying linear operators.

Moreover, not all binding interactions inside the cell occur between two freely diffusing species; binding can also occur between a freely diffusing species, such as RNA polymerase (RNAP), and a spatially fixed species, such as DNA. It has been shown that in this case, spatial heterogeneity modulates the binding affinity between the species [6]. Consequently, it is of interest to extend these results to include fixed species, where the diffusion operator is identically zero.

Finally, in general contraction theory for ODEs, the contraction metric can depend on the state variable. While our results can be easily extended to have $M_1$ depend on $\bar{w}$ by considering virtual displacements of $\bar{w}$ [20], this extension is not straightforward for $M_2$, since dependence on $z^\perp$ breaks the self-adjointness of the differential operators. Future work should aim to allow $M_2$ to depend on $z^\perp$. One potential approach is to project the $z^\perp$ dynamics onto the infinite-dimensional basis of the linear differential operator and have the metric depend on the modal coordinates. However, it will be a challenge to find contraction metrics in that infinite dimensional space.

# 7 Supplementary Information

## 7.1 Available Volume Profiles $v_r(x)$

The diffusion of biomolecules in prokaryotic cells is restricted by the presence of a dense DNA mesh, known as the nucleoid, which contains the cell's genetic material. This phenomenon is quantified using the concept of the available volume profile [2]. Specifically, for a biomolecule with a radius of gyration $r$, the available volume $v_r(x)$ at a spatial location $x \in \Omega$ is given by

$$v_r(x) = e^{-r^2 \rho(x)}, \quad \rho(x) = \frac{1}{1 + e^{20\sqrt{(x-x^*)^T(x-x^*)}}} \tag{55}$$

where $\rho(x)$ represents the density of the nucleoid at location $x$, and $x^*$ denotes the spatial location corresponding to the midpoint of the cell.

## 7.2 Simulation Details For Fig 1 In The Main Text

The numerical method used to discretize and simulate the PDEs is identical to that described in [6]. The parameters for the simulation are $\Omega = (0, 1)$ and 500 spatial nodes. The parameter $\epsilon$ was set to $\epsilon = 1 \times 10^{-2}$. The initial conditions were a uniform profile at unity, specifically $z(0, x) = 1$ for all $x \in (0, 1)$.

For $\epsilon$, the amount of spatial heterogeneity, controlled by $\omega$, was varied. The slope of $\log \|z(t, x)\|$ was calculated after an initial transient, for $80 \leq t \leq 100$. This served as a measure of stability: a positive slope implies instability, whereas a negative slope indicates stability. The critical value of $\omega$ where the slope is zero is denotes in the figure by a vertical dashed red line. The code to replicate this simulation is accessible in: https://github.com/carlobarCodes/PDE_Codes.

## 7.3 Simulation Details for Fig 2 In The Main Text

The numerical method used to discretize and simulate the PDEs is identical to that described in [6]. The parameters for the simulation are $\Omega = (0, 1)$ and 500 spatial nodes. The parameter $x^*$ in (55) is set to $x^* = 1/2$. The initial conditions are $z_1(0, x) = x$ and $z_2(0, x) = 1 + x$.

For Fig. 2-(b), $r = 0$ and $\zeta$ is varied to control the overall speed of diffusion. For Fig. 2-(c), $r$ is varied between zero and unity to evaluate the impact of spatial heterogeneity in the diffusion dynamics of $z_2(t, x)$. For a fixed value of $r$, $\zeta$ was varied between unity and 10 to determine the critical value $\zeta_{\text{cr}}$ at which the dynamics transition from stable to unstable. Specifically, $\zeta_{\text{cr}}$ is defined as the value of $\zeta$ for which the slope of $\log \|z(t, x)\|$ changes from positive to negative. The blue line in Fig. 2-(c) arises from the upper-bound approximation from (43) in the main text. The code to replicate this simulation is accessible in: https://github.com/carlobarCodes/PDE_Codes.

## 7.4 Proof of Lemma 1 In the Main Text

Let $\lambda^*$ be the first non-zero Neumann eigenvalue for the Laplacian on the domain $\Omega$. Consider the operator $L = -\Delta$, defined for all $z \in H^1(\Omega)$ such that $u$ satisfies the Neumann boundary condition:

$$\nabla z \cdot \nu = 0, \quad \forall x \in \partial\Omega,$$

where $\nu$ denotes the outward unit normal vector on $\partial\Omega$. The eigenvalue problem is given by:

$$Lz(x) = -\Delta z(x) = \lambda z(x), \quad \forall x \in \Omega,$$

Then by the min-max principle ([24, Equation 1.37) we have the following inequality

$$\lambda^* \int_\Omega z^\perp(x) z^\perp(x)\, dx \leq \int_\Omega \nabla z^\perp(x) \cdot \nabla z^\perp(x)\, dx, \tag{56}$$

where $z^\perp(x) = z(x) - \int_\Omega z(x) dx$.

Next, for the operator $L(\theta, d, y)$ defined in (13) in the main text, we introduce the coordinate change $u^\perp = \psi_{\theta,d}^{-1}(x) y^\perp(x)$. Under this transformation, the boundary condition becomes:

$$\nabla u^\perp \cdot \nu = 0, \tag{57}$$

where $\nu$ is the outward unit normal vector to the boundary $\partial\Omega$. In this new coordinates, we revisit the inequality from (16) in the main text and apply integrations by parts, the boundary conditions (57), and the inequality (56),

$$\begin{aligned}
\int_\Omega y^\perp(x)\psi_{\theta,d}^{-1}(x)L(\theta,d,y^\perp(x))\,dx &= \int_\Omega u^\perp(x)\mathrm{div}\left(\frac{d^{2\theta}(x)}{\int_\Omega d^{2\theta-1}(x)\,dx}\nabla u^\perp(x)\right)dx \\
&= \int_{\partial\Omega} u^\perp(x)\frac{d^{2\theta}(x)}{\int_\Omega d^{2\theta-1}(x)\,dx}\nabla u^\perp(x)\cdot\nu\,dS - \int_\Omega \frac{d^{2\theta}(x)}{\int_\Omega d^{2\theta-1}(x)\,dx}\nabla u^\perp(x)\cdot\nabla u^\perp(x)\,dx \\
&= -\int_\Omega \frac{d^{2\theta}(x)}{\int_\Omega d^{2\theta-1}(x)\,dx}\nabla u^\perp(x)\cdot\nabla u^\perp(x)\,dx \qquad (58) \\
&\leq -\frac{\min_{x\in\overline{\Omega}} d^{2\theta}(x)}{\int_\Omega d^{2\theta-1}(x)\,dx}\int_\Omega \nabla u^\perp(x)\cdot\nabla u^\perp(x)\,dx \\
&\leq -\lambda^* \frac{\min_{x\in\overline{\Omega}} d^{2\theta}(x)}{\int_\Omega d^{2\theta-1}(x)\,dx}\int_\Omega u^\perp(x)u^\perp(x)\,dx \\
&= -\lambda^* \frac{\min_{x\in\overline{\Omega}} d^{2\theta}(x)}{\int_\Omega d^{2\theta-1}(x)\,dx}\int_\Omega \psi_{\theta,d}^{-1}(x)y^\perp(x)\psi_{\theta,d}^{-1}(x)y^\perp(x)\,dx \\
&\leq -\lambda^* \frac{\min_{x\in\overline{\Omega}} d^{2\theta}(x)}{\int_\Omega d^{2\theta-1}(x)\,dx}\left(\min_{x\in\overline{\Omega}}\psi_{\theta,d}^{-1}(x)\right)\int_\Omega y^\perp(x)\psi_{\theta,d}^{-1}(x)y^\perp(x)\,dx \\
&= -\lambda^* \frac{\min_{x\in\overline{\Omega}} d^{2\theta}(x)}{\max_{x\in\overline{\Omega}} d^{2\theta-1}(x)}\int_\Omega y^\perp(x)\psi_{\theta,d}^{-1}(x)y^\perp(x)\,dx \\
&= -\lambda_{\theta,d}\int_\Omega y^\perp(x)\psi_{\theta,d}^{-1}(x)y^\perp(x)\,dx,
\end{aligned}$$

where

$$\lambda_{\theta,d} = \lambda^* \frac{\min_{x\in\overline{\Omega}} d^{2\theta}(x)}{\max_{x\in\overline{\Omega}} d^{2\theta-1}(x)},$$

and $dS$ denotes the surface measure on the boundary $\partial\Omega$.

### 7.5 Proof Theorem 1 In the Main Text

#### 7.5.1 Preliminary results

**Claim** 1. *Let $y^\perp(t,x) \in \mathbb{R}^m$ be such that*

$$\int_\Omega y^\perp(t,x)\,dx = 0,$$

*and $M_2(t,x)$, $\mathcal{L}_\Theta(D,y^\perp)$, and $\Lambda_{\Theta,D}$ be defined as in the statement of the theorem, then*

$$\int_\Omega y^{\perp,T}M_2\mathcal{L}_\Theta(D,y^\perp)dx \leq \int_\Omega y^{\perp,T}M_2\Lambda_{\Theta,D}y^\perp dx$$

*Proof.* The proof utilizes the digonal nature of $M_2(t,x)$ and the innequlity from Lemma 1 in the main text.

$$\begin{aligned}
\int_\Omega y^{\perp,T}M_2\mathcal{L}_\Theta(D,y^\perp)dx &= \sum_{i=1}^m \Gamma^{i,i}(t)\int_\Omega y^{\perp,i}\psi_{\Theta^i,D^i}(x)L(\Theta^i,D^i,y^{\perp,i})dx \\
&\leq \sum_{i=1}^m \Gamma^{i,i}(t)\int_\Omega y^{\perp,i}\lambda_{\Theta^i,D^i}\psi_{\Theta^i,D^i}(x)y^{\perp,i}dx \\
&= \int_\Omega y^{\perp,T}M_2\Lambda_{\Theta,D}y^\perp dx
\end{aligned}$$

□

**Claim 2.** *Let $h : \mathbb{R}^p \to \mathbb{R}^m$ be a continuously differentiable function with respect to $\xi \in \mathbb{R}^p$. Then, for any $\xi_1, \xi_2 \in \mathbb{R}^p$, the following holds:*

$$h(\xi_2) - h(\xi_1) = \left[\int_0^1 \frac{\partial h}{\partial \xi}(\xi_1 + s(\xi_2 - \xi_1))\, ds\right](\xi_2 - \xi_1).$$

*Proof.* This holds by the mean-value theorem (equation 10 in chapter 3 of [32]). □

**Lemma 2.** *Let $y^\perp(t, x) \in \mathbb{R}^m$ be such that*

$$\int_\Omega y^\perp(t, x)\, dx = 0.$$

*Let $\bar{u}(t) \in \mathbb{R}^n$ be a vector independent of $x$, and let*

$$v(s) \in \chi_v \quad (\text{for } s \in (0, 1)), \quad \text{where } \chi_v \subset \mathbb{R}^{n+m}.$$

*We define the matrix $B(t, x, v(s)) \in \mathbb{R}^{m \times n}$ and let*

$$B^\perp(t, x, v(s)) := B(t, x, v(s)) - \int_\Omega B(t, x, v(s))\, dx.$$

*Assume there exists $\beta \geq 0$ such that for all $t \geq 0$,*

$$\int_\Omega \left[\sup_{v \in \chi_v} B^{\perp, \top}(t, x, v(s))\, B^\perp(t, x, v(s))\right] dx \ \leq\ \beta^2 I_{n \times n}.$$

*Then,*

$$\int_\Omega \int_0^1 \left[y^{\perp,\top}(t,x)\, B(t, x, v(s))\, \bar{u}(t)\right] ds\, dx \ \leq\ \beta \left[\int_\Omega y^{\perp,\top}(t, x)\, y^\perp(t, x)\, dx\right]^{1/2} \left[\bar{u}(t)^\top \bar{u}(t)\right]^{1/2}.$$

*Proof.* First we leverage the fact that

$$B(t, x, v(s)) \ =\ B^\perp(t, x, v(s)) + \int_\Omega B(t, y, v(s))\, dy.$$

Thus,

$$\int_\Omega \int_0^1 \left[y^{\perp,\top}(t,x)\, B(t, x, v(s))\, \bar{u}(t)\right] ds\, dx$$

$$= \int_\Omega \int_0^1 y^{\perp,\top}(t,x) \left[B^\perp(t, x, v(s)) + \int_\Omega B(t, y, v(s))\, dy\right] \bar{u}(t)\, ds\, dx$$

$$= \int_\Omega \int_0^1 \left[y^{\perp,\top}(t,x)\, B^\perp(t, x, v(s))\, \bar{u}(t)\right] ds\, dx + \left(\int_\Omega y^{\perp,\top}(t, x)\, dx\right)\left(\int_0^1 \int_\Omega B(t, y, v(s))\, dy\, ds\right) \bar{u}(t).$$

However,

$$\int_\Omega y^\perp(t, x)\, dx = 0,$$

so the second term vanishes. Therefore,

$$\int_\Omega \int_0^1 \left[y^{\perp,\top}(t,x)\, B(t, x, v(s))\, \bar{u}(t)\right] ds\, dx \ =\ \int_\Omega \int_0^1 \left[y^{\perp,\top}(t,x)\, B^\perp(t, x, v(s))\, \bar{u}(t)\right] ds\, dx. \tag{59}$$

By the usual Cauchy–Schwarz (for vectors in $\mathbb{R}^m$),

$$\left|y^{\perp,\top}(t,x)\left[B^\perp(\cdots)\, \bar{u}(t)\right]\right| \ \leq\ \sqrt{y^{\perp,\top}(t,x)\, y^\perp(t,x)}\, \sqrt{\left[B^\perp(\cdots)\, \bar{u}(t)\right]^\top \left[B^\perp(\cdots)\, \bar{u}(t)\right]}.$$

Then, we can simplify (59) as

$$\int_0^1 \int_\Omega \left| y^{\perp,\top}(t,x) B^\perp(t,x,v(s)) \bar u(t) \right| dx\, ds$$
$$\leq \int_0^1 \int_\Omega \sqrt{y^{\perp,\top}(t,x) y^\perp(t,x)} \sqrt{\left[B^\perp(\cdots)\bar u(t)\right]^\top \left[B^\perp(\cdots)\bar u(t)\right]}\, dx\, ds,$$

where Tonelli/Fubini's theorem was used to swap the order of integrals (the integrand is nonnegative). Next, apply the Cauchy–Schwarz inequality in the $L^2(\Omega)$-space to the integral over $\Omega$ for each fixed $s$:

$$\int_\Omega \sqrt{y^{\perp,\top}(t,x) y^\perp(t,x)} \sqrt{\left[B^\perp(\cdots)\bar u(t)\right]^\top \left[B^\perp(\cdots)\bar u(t)\right]}\, dx$$
$$\leq \sqrt{\int_\Omega y^{\perp,\top}(t,x) y^\perp(t,x)\, dx} \sqrt{\int_\Omega \left[\left[B^\perp(\cdots)\bar u(t)\right]^\top \left[B^\perp(\cdots)\bar u(t)\right]\right] dx}.$$

Note that
$$\left[B^\perp(\cdots)\bar u(t)\right]^\top \left[B^\perp(\cdots)\bar u(t)\right] = \bar u(t)^\top \left[B^{\perp,\top}(t,x,v(s)) B^\perp(t,x,v(s))\right] \bar u(t),$$

so
$$\int_\Omega \left[B^\perp(\cdots)\bar u(t)\right]^\top \left[B^\perp(\cdots)\bar u(t)\right] dx = \bar u(t)^\top \left[\int_\Omega B^{\perp,\top}(t,x,v(s)) B^\perp(t,x,v(s))\, dx\right] \bar u(t).$$

By assumption, for every $t \geq 0$,

$$\int_\Omega \sup_{v \in \chi_v} B^{\perp,\top}(t,x,v) B^\perp(t,x,v)\, dx \leq \beta^2 I_{n \times n}.$$

Hence, for each $v(s) \in \chi_v$,

$$\int_\Omega B^{\perp,\top}(t,x,v(s)) B^\perp(t,x,v(s))\, dx \leq \beta^2 I_{n \times n}.$$

Therefore,
$$\bar u(t)^\top \left(\int_\Omega B^{\perp,\top}(t,x,v(s)) B^\perp(t,x,v(s))\, dx\right) \bar u(t) \leq \beta^2\, \bar u(t)^\top \bar u(t).$$

Combining everything, we obtain

$$\int_0^1 \int_\Omega \left| y^{\perp,\top}(t,x) B^\perp(t,x,v(s)) \bar u(t) \right| dx\, ds$$
$$\leq \int_0^1 \sqrt{\int_\Omega y^{\perp,\top}(t,x) y^\perp(t,x)\, dx} \sqrt{\bar u(t)^\top \left[\int_\Omega B^{\perp,\top}(t,x,v(s)) B^\perp(t,x,v(s))\, dx\right] \bar u(t)}\, ds$$
$$\leq \int_0^1 \sqrt{\int_\Omega y^{\perp,\top}(t,x) y^\perp(t,x)\, dx} \sqrt{\beta^2 \left(\bar u(t)^\top \bar u(t)\right)}\, ds$$
$$= \int_0^1 \beta \sqrt{\int_\Omega y^{\perp,\top}(t,x) y^\perp(t,x)\, dx} \sqrt{\bar u(t)^\top \bar u(t)}\, ds.$$

Since the integrand does not depend on $s$ except for $v(s)$ (which is bounded by the same $\beta^2$-hypothesis for each $s$), we finally get

$$\int_0^1 \int_\Omega y^{\perp,\top}(t,x) B^\perp(t,x,v(s)) \bar u(t)\, dx\, ds \leq \beta \sqrt{\int_\Omega y^{\perp,\top}(t,x) y^\perp(t,x)\, dx} \sqrt{\bar u(t)^\top \bar u(t)}.$$

Recalling the initial decomposition, we conclude

$$\int_\Omega \int_0^1 \left[y^{\perp,\top}(t,x) B(t,x,v(s)) \bar u(t)\right] ds\, dx \leq \beta \left[\int_\Omega y^{\perp,\top}(t,x) y^\perp(t,x)\, dx\right]^{1/2} \left[\bar u(t)^\top \bar u(t)\right]^{1/2}.$$

This completes the proof. □

### 7.5.2 Main Result

Let $[\bar{w}_1, z_1^\perp]^\top \in \bar{\chi}_w \times \chi^\perp$ and $[\bar{w}_2, z_2^\perp]^\top \in \bar{\chi}_w \times \chi^\perp$ be any two solutions to (18), and define the error variables as $\bar{e} = \bar{w}_2 - \bar{w}_1$ and $e^\perp = z_2^\perp - z_1^\perp$. The dynamics of the error variables are given by:

$$\frac{d\bar{e}}{dt} = f_1(t, \bar{w}_2) - f_1(t, \bar{w}_1) + \int_\Omega [g_1(t, x, z_2^\perp) - g_1(t, x, z_1^\perp)] dx,$$

$$= \left[ \int_0^1 \frac{\partial f_1}{\partial \bar{w}} (\bar{w}_1 + s(\bar{w}_2 - \bar{w}_1)) \, ds \right] \bar{e} + \left[ \int_\Omega \int_0^1 [\frac{\partial g_1}{\partial z^\perp} (z_1^\perp + s(z_2^\perp - z_1^\perp))] ds dx \right] e^\perp,$$

and

$$\frac{\partial e^\perp}{\partial t} = \mathcal{L}_\Theta(D, z_2^\perp) - \mathcal{L}_\Theta(D, z_1^\perp) + f_2^\perp(t, x, z_2^\perp) - f_2^\perp(t, x, z_1^\perp) + g_2^\perp(t, x, \bar{w}_2) - g_2^\perp(t, x, \bar{w}_1),$$

$$= \mathcal{L}_\Theta(D, e^\perp) + \left[ \int_0^1 \frac{\partial f_2}{\partial z^\perp} (z_1^\perp + s(z_2^\perp - z_1^\perp)) ds \right] e^\perp + \left[ \int_0^1 \frac{\partial g_2}{\partial \bar{w}} (\bar{w}_1 + s(\bar{w}_2 - \bar{w}_1)) \, ds \right] \bar{e},$$

where we've leveraged the result from Claim 2 to express differences in terms of the Jacobian.

Let the error metric be defined as

$$v(t) = v_1(t) + v_2(t), \text{ where } v_1(t) = \bar{e}^\top M_1 \bar{e} \text{ and } v_2(t) = \int_\Omega \left( e^{\perp,\top} M_2 e^\perp \right) dx.$$

The error metric, $v(t)$, represents a measure of the difference between the two solutions. Notice that we have the following bounds

$$m_{1,*} \bar{e}^\top \bar{e} \leq v_1(t) \leq m_1^* \bar{e}^\top \bar{e},$$

$$m_{2,*} \int_\Omega \left( e^{\perp,\top} e^\perp \right) dx \leq v_2(t) \leq m_2^* \int_\Omega \left( e^{\perp,\top} e^\perp \right) dx,$$

hence $v(t) \to 0$ implies the norm of the error variables also vanishes at the same rate.

The time derivative of the error metric is given by

$$\frac{dv_1(t)}{dt} = 2\bar{e}^\top M_1 \frac{d\bar{e}}{dt} + \bar{e}^\top \dot{M}_1 \bar{e}$$

$$= 2 \left[ \bar{e}^\top M_1 \left[ \int_0^1 \frac{\partial f_1}{\partial \bar{w}} (\bar{w}_1 + s\bar{e}) ds \right] \bar{e} + \bar{e}^\top M_1 \left[ \int_\Omega \int_0^1 [\frac{\partial g_1}{\partial z^\perp} (z_1^\perp + s e^\perp)] e^\perp ds dx \right] \right] + \bar{e}^\top \dot{M}_1 \bar{e}$$

$$= 2 \left[ \int_0^1 \left( \bar{e}^\top \left[ \frac{1}{2} [\dot{M}_1 + M_1 \frac{\partial f_1}{\partial \bar{w}} (\bar{w}_1 + s\bar{e}) + \frac{\partial f_1}{\partial \bar{w}}^\top (\bar{w}_1 + s\bar{e}) M_1 ] \right] \bar{e} \right) ds \right.$$

$$+ \int_\Omega \int_0^1 \left( \bar{e}^\top M_1 [\frac{\partial g_1}{\partial z^\perp} (z_1^\perp + s e^\perp)] e^\perp \right) ds dx \right]$$

$$\leq 2 \left[ \int_0^1 \left( -\lambda_1 \bar{e}^\top M_1 \bar{e} \right) ds + \int_\Omega \int_0^1 \left( \bar{e}^\top M_1 [\frac{\partial g_1}{\partial z^\perp} (z_1^\perp + s e^\perp)] e^\perp \right) ds dx \right]$$

$$\leq -2\lambda_1 v_1(t) + 2 \int_\Omega \int_0^1 \left( \bar{e}^\top M_1 [\frac{\partial g_1}{\partial z^\perp} (z_1^\perp + s e^\perp)] e^\perp \right) ds dx,$$

where we've leveraged assumption (19) in the theorem statement and the convexity of $\bar{\chi}_w$. Similarly,

$$\begin{aligned}\frac{dv_2(t)}{dt} &= \int_\Omega \left(2e^{\perp,\top} M_2 \frac{\partial e^\perp}{\partial t} + e^{\perp,\top} \dot{M}_2 e^\perp\right) dx \\ &= 2\int_\Omega \left(e^{\perp,\top} M_2 \mathcal{L}_\Theta(D, e^\perp) + e^{\perp,\top} M_2 \left[\int_0^1 \frac{\partial f_2}{\partial z^\perp}(z_1^\perp + se^\perp) ds\right] e^\perp + \frac{1}{2} e^{\perp,\top} \dot{M}_2 e^\perp \right. \\ &\quad \left. + e^{\perp,\top} M_2 \left[\int_0^1 \frac{\partial g_2}{\partial \bar{w}}(\bar{w}_1 + s\bar{e}) ds\right] \bar{e}\right) dx \\ &\leq 2\int_\Omega \int_0^1 \left(e^{\perp,\top} \frac{1}{2}\left[\dot{M}_2 + M_2\left[\frac{\partial f_2}{\partial z^\perp}(z_1^\perp + se^\perp) - \Lambda_{\Theta,D}\right] + \left[\frac{\partial f_2^T}{\partial z^\perp}(z_1^\perp + se^\perp) M_2 - \Lambda_{\Theta,D}\right]\right] e^\perp \\ &\quad + e^{\perp,\top} M_2 \left[\frac{\partial g_2}{\partial \bar{w}}(\bar{w}_1 + s\bar{e})\right] \bar{e}\right) dsdx \\ &\leq 2\int_\Omega \int_0^1 \left(-\lambda_2 e^{\perp,\top} M_2 e^\perp + e^{\perp,\top} M_2 \left[\frac{\partial g_2}{\partial \bar{w}}(\bar{w}_1 + s\bar{e})\right] \bar{e}\right) dsdx \\ &\leq -2\lambda_2 v_2(t) + 2\int_\Omega \int_0^1 \left(e^{\perp,\top} M_2 \left[\frac{\partial g_2}{\partial \bar{w}}(\bar{w}_1 + s\bar{e})\right] \bar{e}\right) dsdx,\end{aligned}$$

where we've leveraged the result from Claim 1, assumption (20) in the theorem statement, and the convexity of $\chi^\perp$. Altogether, this implies that

$$\begin{aligned}\frac{dv(t)}{dt} &= \frac{dv_1(t)}{dt} + \frac{dv_2(t)}{dt} \\ &\leq -2\big(\lambda_1 v_1(t) + \lambda_2 v_2(t)\big) + 2\int_\Omega \int_0^1 \left(e^{\perp,\top} \left[\frac{\partial g_1^T}{\partial z^\perp}(z_1^\perp + se^\perp) M_1 + M_2 \frac{\partial g_2}{\partial \bar{w}}(\bar{w}_1 + s\bar{e})\right] \bar{e}\right) dsdx.\end{aligned}$$

To simplify this expression we apply result from Lemma 2 with $B = G = \left(\frac{\partial g_1}{\partial z^\perp}\right)^\top M_1 + M_2 \frac{\partial g_2}{\partial \bar{w}}$, assumption (21), and the convexity of $\bar{\chi}_w$ and $\chi^\perp$, more precisely,

$$\begin{aligned}\frac{dv(t)}{dt} &\leq -2\big(\lambda_1 v_1(t) + \lambda_2 v_2(t)\big) + 2\beta \left[\int_\Omega \left(e^{\perp,\top} e^\perp\right) dx\right]^{1/2} \left[\bar{e}^\top \bar{e}\right]^{1/2} \\ &\leq -2\big(\lambda_1 v_1(t) + \lambda_2 v_2(t)\big) + 2\beta \left[\frac{v_1(t)}{m_{1,*}}\right]^{1/2} \left[\frac{v_2(t)}{m_{2,*}}\right]^{1/2} \\ &= -2 \begin{bmatrix} v_1^{1/2}(t) & v_2^{1/2}(t) \end{bmatrix} \begin{bmatrix} \lambda_1 & -\sigma \\ -\sigma & \lambda_2 \end{bmatrix} \begin{bmatrix} v_1^{1/2}(t) \\ v_2^{1/2}(t) \end{bmatrix} \\ &\leq -2\lambda^* \big[v_1(t) + v_2(t)\big] \\ &= -2\lambda^* v(t),\end{aligned}$$

where $\sigma = \frac{\beta}{2\sqrt{m_{1,*} m_{2,*}}}$ and $\lambda^*$ is the smallest eigenvalue of $\begin{bmatrix} \lambda_1 & -\sigma \\ -\sigma & \lambda_2 \end{bmatrix}$ and is given by

$$\lambda^* := \left[\frac{\lambda_1 + \lambda_2}{2} - \sqrt{\left(\frac{\lambda_1 - \lambda_2}{2}\right)^2 + \sigma^2}\right].$$

The positivity of $\lambda^*$ requires

$$\lambda_1 \lambda_2 > \sigma^2,$$

and hence the motivation for assumption (22). Thus, altogether this implies that $v(t) \to 0$ exponentially fast since

$$v(t) \leq v(0) e^{-2\lambda^*(t)}.$$

This concludes the proof.

## 7.6 Details for Example 5.5 In the Main Text

### 7.6.1 Positiveness and Boundedness of Concentrations

Intuitively, one would expect intracellular concentrations to be positive and bounded. To this end, we apply the results of [33] to (44) in the main text to derive pointwise bounds on the concentrations of the system and to demonstrate positivity. In Theorem 1 of [33], the authors demonstrated that if the diffusion term is described by a general linear elliptic operator in standard form, then a convex set $S \subset \mathbb{R}^n$ is positively invariant if the reaction vector field points inward at the boundary of $S$. This is analogous to the case where the dynamics are described by ODEs, in which showing that the dynamic vector points inward at the boundaries of the set ensures invariance.

However, the system (44) in the main text is not in the form specified in [33] (Equation 1.1), since the diffusion operator is not in standard elliptic form. More precisely,

$$L(\theta, d, y) = d(x)\Delta y(x) + 2(1-\theta)\nabla d(x) \cdot \nabla y(x) + (1-2\theta)y(x)\Delta d(x),$$

and hence the sign-indefinite term $(1-2\theta)y(x)\Delta d(x)$ prevents the application of the invariance theorem.

To transform the operator into standard elliptic form, we define the change of variables:

$$u_m(\tau, x) := \frac{m(\tau, x)}{v_{r_m}(x)}, \quad u_r(\tau, x) := \frac{R(\tau, x)}{v_{r_r}(x)}, \quad u_c(\tau, x) := \frac{c_r(\tau, x)}{v_{r_c}(x)}.$$

Under this coordinate change, the dynamics are given by:

$$\begin{aligned}
\frac{\partial u_m}{\partial \tau} &= \chi_m \big[ v_{r_m} \Delta u_m + 2\nabla v_{r_m} \cdot \nabla u_m \big] - v_{r_r}\Big(\frac{u_m u_r}{K} - u_c\Big), \\
\frac{\partial u_r}{\partial \tau} &= \chi_r \big[ v_{r_r} \Delta u_r + 2\nabla v_{r_r} \cdot \nabla u_r \big] - v_{r_m}\Big(\frac{u_m u_r}{K} - u_c\Big), \\
\frac{\partial u_c}{\partial \tau} &= \chi_c \big[ v_{r_c} \Delta u_c + 2\nabla v_{r_c} \cdot \nabla u_c \big] + \Big(\frac{u_m u_r}{K} - u_c\Big),
\end{aligned} \quad (60)$$

where we used the fact that $v_{r_c}(x) = v_{r_r}(x)v_{r_m}(x)$, which follows from the factorization of the polysome available volume into independent contributions from ribosome-free mRNA and ribosomes, as shown in Equation (S12) of [2]). The boundary conditions under this coodinate change are of Neumann type:

$$\nabla u_m \cdot \nu = 0, \quad \nabla u_r \cdot \nu = 0, \quad \nabla u_c \cdot \nu = 0, \quad \forall x \in \partial\Omega, \, \forall t > 0.$$

The diffusion operator in (60) is in standard elliptic form since it only contains the Laplacian and gradient of the state variable.

Next, for any finite $C^* > 0$, we define the following convex set $S$ and demonstrate that it is invariant:

$$S = \big\{ (u_m, u_r, u_c) \in \mathbb{R}^3 \,\big|\, u_m \geq 0,\, u_r \geq 0,\, u_c \geq 0,\, v_{r_m}(x)\,u_m + v_{r_r}(x)\,u_r + 2\,v_{r_c}(x)\,u_c \leq C^* \big\}.$$

We define the reaction-field vector $\mathbf{F}$ as:

$$\mathbf{F} = \Big( -v_{r_r}\big(\tfrac{u_m u_r}{K} - u_c\big),\, -v_{r_m}\big(\tfrac{u_m u_r}{K} - u_c\big),\, \tfrac{u_m u_r}{K} - u_c \Big).$$

The boundary of $S$ is composed of the following components, along with their respective normal vectors and inner products:

1. Plane: $u_m = 0$ (for all $u_r, u_c \geq 0$),
   Normal vector: $\mathbf{n} = (-1, 0, 0)$,
   Inner product with reaction field: $\mathbf{n} \cdot \mathbf{F} = v_r \left( \frac{u_m u_r}{K} - u_c \right) = -v_r u_c \leq 0$.

2. Plane: $u_r = 0$ (for all $u_m, u_c \geq 0$),
   Normal vector: $\mathbf{n} = (0, -1, 0)$,
   Inner product with reaction field: $\mathbf{n} \cdot \mathbf{F} = v_{r_m} \left( \frac{u_m u_r}{K} - u_c \right) = -v_{r_m} u_c \leq 0$.

3. Plane: $u_c = 0$ (for all $u_m, u_r \geq 0$),
   Normal vector: $\mathbf{n} = (0, 0, -1)$,
   Inner product with reaction field: $\mathbf{n} \cdot \mathbf{F} = -\left( \frac{u_m u_r}{K} - u_c \right) = -\frac{u_m u_r}{K} \leq 0$.

4. Plane: $v_{r_m}(x) u_m + v_{r_r}(x) u_r + 2 v_{r_c}(x) u_c = C^*$,
   Normal vector: $\mathbf{n} = \left( v_{r_m}(x), v_{r_r}(x), 2 v_{r_c}(x) \right)$,
   Inner product with reaction field: $\mathbf{n} \cdot \mathbf{F} = 0$.

The normal vector for plane 4, $v_{r_m}(x) u_m + v_{r_r}(x) u_r + 2 v_{r_c}(x) u_c = C^*$, depends on $x$, which is not consistent with the statement of Theorem 1 in [33]. However, motivated by Theorem 5 in their work, we can simply enlarge the state space so that $v_{r_m}(x)$, $v_{r_r}(x)$, and $v_{r_c}(x)$ become additional dependent variables:

$$u_4(\tau, x) = v_{r_m}(x), \quad u_5(\tau, x) = v_{r_r}(x), \quad u_6(\tau, x) = v_{r_c}(x).$$

We append the trivial PDEs

$$\frac{\partial u_4}{\partial \tau} = 0, \quad \frac{\partial u_5}{\partial \tau} = 0, \quad \frac{\partial u_6}{\partial \tau} = 0$$

to form a weakly coupled parabolic system in 6 variables. In this extended viewpoint, the set $S$ is described by

$$u_m, u_r, u_c \geq 0 \quad \text{and} \quad u_m u_4 + u_r u_5 + 2 u_c u_6 \leq C^*.$$

One checks (using the normal vectors above) that the reaction/diffusion field is inward-pointing on each boundary facet. Thus, by Weinberger's criterion(specifically Theorem 5 in [33]), $S$ is invariant under the dynamics, and so the original variables $(u_m, u_r, u_c)$ remain in $S$ for all time if initially so.

Let

$$v^*_{r_m} := \sup_{x \in \Omega} v_{r_m}(x), \quad v^*_{r_r} := \sup_{x \in \Omega} v_{r_r}(x), \quad v^*_{r_c} := \sup_{x \in \Omega} v_{r_c}(x),$$

and

$$v_{r_m,*} := \inf x \in \Omega v_{r_m}(x), \quad v_{r_r,*} := \inf_{x \in \Omega} v_{r_r}(x), \quad v_{r_c,*} := \inf_{x \in \Omega} v_{r_c}(x).$$

Altogether this implies that if the initial conditions satisfy

$$\frac{v^*_{r_m}}{v_{r_m,*}} m(0, x) + \frac{v^*_{r_r}}{v_{r_r,*}} R(0, x) + \frac{v^*_{r_c}}{v_{r_c,*}} c_r(0, x) \leq C^*, \quad \forall x \in \Omega$$

then the trajectories of the system are bounded by

$$\frac{v_{r_m,*}}{v^*_{r_m}} m(\tau, x) + \frac{v_{r_r,*}}{v^*_{r_r}} R(\tau, x) + \frac{v_{r_c,*}}{v^*_{r_c}} c_r(\tau, x) \leq C^*, \quad \forall \tau \geq 0, \quad \forall x \in \Omega.$$

We now define bounds on $m(\tau, x)$ and $R(\tau, x)$ as

$$m^* := C^* \frac{v^*_{r_m}}{v_{r_m,*}}, \quad \text{and} \quad R^* := C^* \frac{v^*_{r_r}}{v_{r_r,*}}. \tag{61}$$

From the positivity of the state variables, it implies that the space average variables are bounded by:

$$0 \leq \bar{m}(\tau) \leq m_T, \quad \text{and} \quad 0 \leq \bar{R}(\tau) \leq R_T, \quad \forall t \geq 0, \tag{62}$$

where $m_T$ and $R_T$ are the total space average concentration defined in (46) in the main text. Finally, we have that

$$0 \leq \hat{v}_{r_m} \bar{m}(\tau) + m^{\perp}(\tau, x) \leq C^* \frac{v^*_{r_m}}{v_{r_m,*}} \quad \text{and} \quad 0 \leq \hat{v}_{r_r} \bar{R}(\tau) + R^{\perp}(\tau, x) \leq C^* \frac{v^*_{r_r}}{v_{r_r,*}},$$

this implies that if we let $m^{\perp,*} := \max(\hat{v}^*_{r_m} \bar{m}_T, C^* \frac{v^*_{r_m}}{v_{r_m,*}})$ and $R^{\perp,*} := \max(\hat{v}^*_{r_r} \bar{R}_T, C^* \frac{v^*_{r_r}}{v_{r_r,*}})$

$$|m^{\perp}(\tau, x)| \leq m^{\perp,*} \quad \text{and} \quad |R^{\perp}(\tau, x)| \leq R^{\perp,*}, \quad \forall \tau \geq 0, \quad \forall x \in \Omega. \tag{63}$$

### 7.6.2 Satisfying the Contraction Conditions

Next, we revisit the contraction conditions of Example 5.5 in the main text. But first we demonstrate the boundedness of $\eta(\tau)$ and $\tilde{u}(\tau, x)$. From the positivity of $\bar{m}(\tau)$ and $\bar{R}(\tau)$ we have that

$$1 \leq \eta(\tau) \leq 1 + \theta \frac{(\bar{m}_T + \bar{R}_T)}{K}, \quad \forall t \geq 0, \tag{64}$$

Next, notice that

$$|v_{r_m}(x)\bar{m}(\tau) + \frac{1}{2}m^{\perp}(\tau, x)| \leq \hat{v}^*_{r_m} \bar{m}_T + \frac{1}{2}m^{\perp,*}$$

and

$$|v_{r_r}(x)\bar{R}(\tau) + \frac{1}{2}R^{\perp}(\tau, x)| \leq \hat{v}^*_{r_r} \bar{R}_T + \frac{1}{2}R^{\perp,*}$$

Let $\beta^2_u := \left(\hat{v}^*_{r_m} \bar{m}_T + \frac{1}{2} m^{\perp,*}\right)^2 + \left(\hat{v}^*_{r_r} \bar{R}_T + \frac{1}{2} R^{\perp,*}\right)^2 + 1$ and thus

$$\left(0 \leq \tilde{u}^T(\tau, x)\tilde{u}(\tau, x) \leq \beta^2_u, \quad \text{and} \quad 0 \leq u^T(\tau, x)u(\tau, x) \leq \beta^2_u - 1\right), \quad \forall \tau \geq 0, \quad \forall x \in \Omega,$$

$$v^T v = 3$$

1. *Stability of the Spatial Average Dynamics* (19): We have $\frac{\partial f_1}{\partial \bar{y}_w} = -\eta(\tau)$. Letting $M_1 = 1$, we approximate a bound for $\lambda_1$:

   $$\frac{1}{2}\left[-\eta(\tau) - \eta^T(\tau)\right] \leq -\lambda_1, \quad \forall \tau \geq 0,$$

   which is satisfied with $\lambda_1 = 1$ due to the positivity of $\bar{m}(\tau)$ and $\bar{R}(\tau)$.

2. *Stability of the Spatial Deviation Dynamics* (20): We have $\frac{\partial f_2}{\partial y^{\perp}_z} = v\tilde{u}^T(\tau, x)$ and letting $\Gamma(t) = I_{2,2}$, we next approximate a bound for $\lambda_2$ in

   $$\frac{1}{2}\left[\Psi^{-1}_{\Theta,D}(x)\left(v\tilde{u}^T(\tau, x) - \Lambda_{\Theta,D}\right) + \left(\tilde{u}(\tau, x)v^T - \Lambda_{\Theta,D}\right)\Psi^{-1}_{\Theta,D}(x)\right] \leq -\lambda_2 \Psi^{-1}_{\Theta,D}(x), \quad \forall \tau \geq 0, \forall x \in \Omega.$$

   For any $w \in \mathbb{R}^m$ we have that

   $$w^T \left[\frac{1}{2}\left[\Psi^{-1}_{\Theta,D}(x)\left(v\tilde{u}^T(\tau, x) - \Lambda_{\Theta,D}\right) + \left(\tilde{u}(\tau, x)v^T - \Lambda_{\Theta,D}\right)\Psi^{-1}_{\Theta,D}(x)\right]\right] w$$

   $$= -w^T \Lambda_{\Theta,D} \Psi^{-1}_{\Theta,D}(x) w + w^T \tilde{u}(\tau, x) v^T \Psi^{-1}_{\Theta,D}(x) w$$

   $$\leq -\left(\min_{i=\{1,2,3\}} \Lambda^{i,i}_{\Theta,D}\right) w^T \Psi^{-1}_{\Theta,D}(x) w + \left(\Psi^{-1/2}_{\Theta,D}(x) w\right)^T \left(\Psi^{1/2}_{\Theta,D}(x)\tilde{u}(\tau, x)\right) \left(\Psi^{-1/2}_{\Theta,D}(x) v\right)^T \left(\Psi^{-1/2}_{\Theta,D}(x) w\right)$$

   $$\leq \left[-\Lambda_{\Theta,D*} + ||\Psi^{1/2}_{\Theta,D}(x)\tilde{u}(\tau, x)|| ||\Psi^{-1/2}_{\Theta,D}(x) v||\right] w^T \Psi^{-1}_{\Theta,D}(x) w$$

   $$\leq \left[-\Lambda_{\Theta,D*} + \sqrt{\frac{\max_{i=\{1,2,3\}}(\sup_{x \in \Omega} \Psi^{i,i}_{\Theta,D}(x))}{\min_{i=\{1,2,3\}}(\inf_{x \in \Omega} \Psi^{i,i}_{\Theta,D}(x))}} ||\tilde{u}(\tau, x)|| ||v||\right] w^T \Psi^{-1}_{\Theta,D}(x) w$$

   $$\leq \left[-\Lambda_{\Theta,D*} + \sqrt{3 \frac{\Psi^*_{\Theta,D}}{\Psi_{\Theta,D*}}} \beta_u\right] w^T \Psi^{-1}_{\Theta,D}(x) w$$

Thus an estimate for $\lambda_2$ is given by

$$\lambda_2 = \Lambda_{\Theta,D*} - \sqrt{3\frac{\Psi^*_{\Theta,D}}{\Psi_{\Theta,D*}}}\beta_u$$

3. *Boundedness of the Coupling Terms* (21): We have $\frac{\partial g_1}{\partial y_z^\perp} = \eta(t)u^T(\tau,x)$ and $\frac{\partial g_2}{\partial \bar{y}_w} = -\hat{v}_{r_c}(x)v$ and thus from the coupling term equation (23) we have that

$$G(\tau,x) = \eta(\tau)u(\tau,x) + \psi^{-1}_{\Theta,D}(x)\hat{v}_{r_c}(x)v = \eta(\tau)u(\tau,x) + \begin{bmatrix} \frac{\hat{v}_{r_r}(x)}{\theta} & \frac{\hat{v}_{r_m}(x)}{\theta} & -1 \end{bmatrix}^T$$

$$G^\perp(\tau,x) = u_\eta^\perp(\tau,x) + h(x),$$

where

$$u_\eta^\perp(\tau,x) := \eta(t)\begin{bmatrix} \frac{[\hat{v}_{r_r}(x)-1]\bar{R}(\tau)+\frac{1}{2}R^\perp(\tau,x)}{K}, & \frac{[\hat{v}_{r_m}(x)-1]\bar{m}(\tau)+\frac{1}{2}m^\perp(\tau,x)}{K}, & 0 \end{bmatrix}^T$$

and $h(x) := \begin{bmatrix} \frac{\hat{v}_{r_r}(x)-1}{\theta} & \frac{\hat{v}_{r_m}(x)-1}{\theta} & 0 \end{bmatrix}^T$. The norm of these quantities is bounded by

$$\|u_\eta^\perp(\tau,x)\|^2 = \frac{\eta(\tau)}{K^2}\int_\Omega \left(\left([\hat{v}_{r_r}(x)-1]\bar{R}(\tau) + \frac{1}{2}R^\perp(\tau,x)\right)^2 + \left([\hat{v}_{r_m}(x)-1]\bar{m}(\tau) + \frac{1}{2}m^\perp(\tau,x)\right)^2\right)dx$$

$$\leq \frac{1}{K^2}\left(1 + \theta\frac{(\bar{m}_T + \bar{R}_T)}{K}\right)\int_\Omega \left(\left([\hat{v}_{r_r}(x)-1]\bar{R}_T + \frac{1}{2}R^{\perp,*}\right)^2 + \left([\hat{v}_{r_m}(x)-1]\bar{m}_T + \frac{1}{2}m^{\perp,*}\right)^2\right)dx$$

and we define

$$\beta_{u,\eta}^2 := \frac{1}{K^2}\left(1 + \theta\frac{(\bar{m}_T + \bar{R}_T)}{K}\right)\int_\Omega \left(\left([\hat{v}_{r_r}(x)-1]\bar{R}_T + \frac{1}{2}R^{\perp,*}\right)^2 + \left([\hat{v}_{r_m}(x)-1]\bar{m}_T + \frac{1}{2}m^{\perp,*}\right)^2\right)dx.$$

Similarly,

$$\beta_h^2 := \|h^T(x)\|^2 = \frac{1}{\theta^2}\int_\Omega (\hat{v}_{r_r}(x)-1)^2 + (\hat{v}_{r_m}(x)-1)^2)dx.$$

Therefore,

$$\int_\Omega G^{\perp,T}G^\perp dx = \int_\Omega \left(u_\eta^\perp(\tau,x) + h(x)\right)^T\left(u_\eta^\perp(\tau,x) + h(x)\right)dx$$

$$\leq \int_\Omega \left(\sqrt{u_\eta^{\perp,\top}(\tau,x)u_\eta^\perp(\tau,x)} + \sqrt{h^T(x)h(x)}\right)^2 dx$$

$$\leq \left(\sqrt{\int_\Omega u_\eta^{\perp,\top}(\tau,x)u_\eta^\perp(\tau,x)dx} + \sqrt{\int_\Omega h^T(x)h(x)dx}\right)^2$$

$$= \left(\|u_\eta^\perp(\tau,x)\| + \|h^T(x)\|\right)^2$$

$$\leq \left(\beta_{u,\eta} + \beta_h\right)^2.$$

Thus $\beta = \beta_{u,\eta} + \beta_h$ satisfies the inequality.

4. *Small Gain Condition* (22):

$$\lambda_1\lambda_2 > \frac{\beta^2}{4} \implies \left(\Lambda_{\Theta,D*} - \sqrt{3\frac{\Psi^*_{\Theta,D}}{\Psi_{\Theta,D*}}}\beta_u\right) > \left(\frac{\beta_{u,\eta}+\beta_h}{2}\right)^2 \implies \Lambda_{\Theta,D*} > \left(\frac{\beta_{u,\eta}+\beta_h}{2}\right)^2 + \sqrt{3\frac{\Psi^*_{\Theta,D}}{\Psi_{\Theta,D*}}}\beta_u$$